\documentclass{svjour3}
\usepackage{latexsym}
\usepackage{amssymb}
\usepackage{dsfont}
\usepackage{verbatim}
\usepackage{comment}
\usepackage{tikz}
\usepackage{pgf}
\usepackage{bm}
\usepackage{graphicx,psfrag}
\usepackage{latexsym}
\usepackage{multirow}
\usepackage{subfig}

\newcommand{\kleinsbug}{(3.2)}
\newcommand{\wabs}[1]{\left|#1\right|}

\newcommand{\wfh}[1]{\wfc{r_{\wvec{x},\delta}[\wvec{y}]}{#1}}
\newcommand{\wfhk}[1]{\wfc{r^{(k)}_{\wvec{x},d}[\wvec{y}]}{#1}}
\newcommand{\wfhi}{r_{\wvec{x},\delta}}
\newcommand{\wxfh}[1]{\wfc{\tilde{r}_{\wvec{x},d, \tilde{n},\tilde{d}}[\wvec{y}]}{#1}}
\newcommand{\wxfhdiag}[1]{\wfc{\tilde{r}_{\wvec{x},d,d,d}[\wvec{y}]}{#1}}
\newcommand{\wxfhi}{\tilde{r}_{\wvec{x},d, \tilde{n},\tilde{d}}}

\newcommand{\wlebcfh}{\tilde{\Lambda}_{\wvec{x},\delta}}
\newcommand{\wlebffh}[1]{\wfc{\Lambda_{\wvec{x},\delta}}{#1}}
\newcommand{\wlebcxfh}{\tilde{\Lambda}_{\wvec{x},d,\tilde{n},\tilde{d}}}
\newcommand{\wleby}{\tilde{\Lambda}_{\wvec{x},d,Y}}
\newcommand{\wlebcdiag}{\tilde{\Lambda}_{\wvec{x},d,d,d}}
\newcommand{\wlebfxfh}[1]{\wfc{\tilde{\Lambda}_{\wvec{x},d,\tilde{n},\tilde{d}}}{#1}}
\newcommand{\wlebfxfhdiag}[1]{\wfc{\tilde{\Lambda}_{\wvec{x},d,d,d}}{#1}}

\newcommand{\wrone}{\mathds R}

\newcommand{\wfc}[2]{{#1}\!\left(#2\right)}

\newcommand{\wi}[1]{\wrm{i}}

\newcommand{\wlr}[1]{\left( #1 \right)}

\newcommand{\wnorm}[1]{\left\| {#1} \right\|}

\newcommand{\wref}[1]{(\ref{#1})}
\newcommand{\wrn}[1]{{\mathds R}^{#1}}

\newcommand{\wrm}[1]{\mathrm{#1}}

\newcommand{\wrounde}[1]{\wfc{\wrm{fl}}{#1}}

\newcommand{\wset}[1]{{\left\{ #1 \right\}}}

\newcommand{\wvec}[1]{\mathbf{#1}}
\newcommand{\wvecB}[2]{{\left( \begin{array}{c}  #1 \\ #2 \end{array} \right)}}

\newcommand{\wz}{\mathds Z}

\smartqed
\begin{document}

\title{The stability of extended Floater-Hormann interpolants}
\author{Andr\'{e} Pierro de Camargo and Walter F. Mascarenhas}

\institute{Andr\'{e} Pierro de Camargo at Centro de Matem\'{a}tica, Computa\c{c}\~{a}o e
Cogni\c{c}\~{a}o, Universidade Federal do ABC – UFABC, Rua Santa Ad\'{e}lia, 166, bairro Bangu, CEP 09210-170 Santo Andr\'{e}, SP, Brazil,
and Walter F. Mascarenhas at  Instituto de Matem\'{a}tica e Estat\'{i}stica, Universidade de S\~{a}o Paulo,
           Cidade Universit\'{a}ria, Rua do Mat\~{a}o 1010, S\~{a}o Paulo SP, Brazil. CEP 05508-090
              Tel.: +55-11-3091 5411, Fax: +55-11-3091 6134,   \email{walter.mascarenhas@gmail.com}.
              Andr\'{e} is supported by grant 14225012012-0 from
              Conselho Nacional de Desenvolvimento Cient\'{i}fico e Tecnol\'{o}gico, CNPq.
              Walter is supported by grant 2013/10916-2 from Funda\c{c}\~{a}o de Amparo \`{a}
              Pesquisa do Estado de S\~{a}o Paulo (FAPESP.) }

\maketitle

%    \subjclass is required.
%\subjclass[2010]{Primary 65D05, 41A20, 41A25.}

%\keywords{Rational interpolation, barycentric formula, Lebesgue constant}
%date{}

\begin{abstract}
We present a new analysis of the stability of extended
Floater-Hormann interpolants, in which  both noisy data
and rounding errors are considered. Contrary to what is claimed in the
current literature, we show that the Lebesgue constant of these interpolants
can grow exponentially with the parameters that define them,
and we emphasize the importance of using the proper interpretation of
the Lebesgue constant in order to estimate correctly the effects of noise
and rounding errors.
We also present a simple condition that implies the backward instability of the
barycentric formula used to implement extended interpolants.
Our experiments show that extended interpolants mentioned
in the literature satisfy this condition and, therefore, the formula
used to implement them is not backward stable.
Finally, we explain that the
extrapolation step is a significant source of numerical instability for extended
interpolants based on extrapolation.
\end{abstract}

\maketitle

\section{Introduction}
Given nodes $\wvec{x} = \wlr{x_0,\dots,x_n}$, an integer $\delta$ with  $0 \leq \delta \leq n$ and function values $\wvec{y} = \wlr{y_0,\dots,y_n}$,
the Floater-Hormann interpolation formula is defined as
\begin{equation}
\label{fhdef}
\wfc{r_\delta}{t,\wvec{x},\wvec{y}} := \frac{\sum_{i = 0}^{n-\delta} \wfc{\lambda_i}{t, \wvec{x}} \wfc{p_i}{t, \wvec{x}, \wvec{y}}}
{\sum_{i = 0}^{n-\delta} \wfc{\lambda_i}{t,\wvec{x}}},
\end{equation}
where $\wfc{p_i}{t,\wvec{x},\wvec{y}}$ is the unique polynomial of degree
at most $\delta$ which interpolates $y_i, y_{i+1}, \dots, y_{i+\delta}$ at $x_i, x_{i+1}, \dots, x_{i+\delta}$, and the weights $\lambda_i$ are
defined as
\[
\wfc{\lambda_i}{t,\wvec{x}} := \frac{\wlr{-1}^i}{\wlr{t - x_i} \wlr{t - x_{i+1}} \dots \wlr{t - x_{i+\delta}}}, \hspace{1cm} \wrm{for} \ i = 0, \dots, n - \delta.
\]

In exact arithmetic, when $y_i = \wfc{f}{x_i}$ for a smooth function $f$, the
error incurred by the Floater-Hormann interpolant defined by $\wvec{x}$
and $\delta$ is of order $h^{\delta+1}$, where
\[
h := \max_{0 \leq k < n}  x_{k+1} - x_k.
\]
Unfortunately, when the nodes are equally spaced the Lebesgue constant
of the Floater-Hormann interpolant defined by $\wvec{x}$ and $\delta$ grows exponentially with $\delta$
(see \cite{Bos}).
Therefore, $\delta$ must be chosen carefully in order to balance the high order of approximation $h^{\delta+1}$
with the numerical errors due to large Lebesgue constants.

In an attempt to reduce the effects of the large Lebesgue constants for equally
spaced nodes, Klein \cite{Klein} introduced the {\it extended Floater-Hormann} interpolants.
These interpolants are defined in terms of an integer parameter $d \geq 0$,
extended nodes
$\tilde{x}_{-d}, \tilde{x}_{1 - d}, \dots, \tilde{x}_0, \dots \tilde{x}_n, \dots, \tilde{x}_{n + d}$,
with $\tilde{x}_i = x_0 + i h$, and
extended function values $\tilde{y}_{-d}, \dots \tilde{y}_0, \dots \tilde{y}_n, \dots \tilde{y}_{n + d}$.
Each
\[
\tilde{\wvec{y}} = \wfc{Y}{d,\wvec{x},\wvec{y}} \in \wrn{n + 2 d + 1}
\]
in combination with $r_d$ in \wref{fhdef} leads to an extended interpolant given by
\begin{equation}
\label{defrdy}
\wfc{\tilde{r}_{d, Y}}{t, \wvec{x}, \wvec{y}} :=
 \wfc{r_d}{t, \tilde{\wvec{x}}, \wfc{Y}{d,\wvec{x},\wvec{y}}}.
\end{equation}
The choice of the extended function values is a crucial point regarding
the stability and accuracy of the extended interpolants. Usually, we do not have
information outside of the interpolation interval
and in practice the $\tilde{y}$ must be estimated, and they will not be exact.
To the best of our knowledge, the only concrete way for choosing the $\tilde{y}$ mentioned in the literature prior to our writing
of this article is the one outlined in the fifth page of Klein and Berrut \cite{BerrutKleinCAM},
which is based on two additional parameters $\tilde{d}$ and $\tilde{n}$:
\begin{quote}
``More precisely, $2d$ extra nodes $x_{-d},\dots,x_{-1}$
are considered, $d$ on each side of the interval, and approximate values $\tilde{f}_i$ of  $f$ at these nodes are computed by a discrete
Taylor polynomial with derivatives approximated by (linear rational) finite differences (see Section 8) {\it using only the given values of
$f$ in $[a,b]$.} These finite differences are the derivatives of the Floater-Hormann family
with parameters $\tilde{d}$ in the nodes $x_0,\dots,x_{\tilde{n}}$, resp. $x_{n-\tilde{n}}, \dots, x_n$,
for an $\tilde{n}$ much smaller than $n$. At the original nodes, $x_j$, $j = 0,\dots,n$, the given $f_f$ are used.''
\end{quote}
Klein \cite{Klein} shows that the order of approximation of the
extended interpolant above is $h^{\mu + 1}$,
where $\mu := \min\{d, \tilde{d}\}$. Usual Floater-Hormann interpolants
have order of approximation $\delta$, and
$\mu$ is the analogous to the parameter $\delta$ used to define usual
Floater-Hormann interpolants. Therefore, it is important to distinguish $\delta$ from $d$.
In fact, when choosing the parameters in practice, one must
be aware that the order of approximation $\mu$ will be unaffected
by increasing the parameter $d$ once this parameter is already larger
than $\tilde{d}$. This argument and our practical experience with
extended interpolants suggest that as first choice
one should pick $d = \tilde{d} = \delta$ (see also Fig. \ref{figure_conjecture} below.)
For this reasons, our theory pays special attention to the case $d = \delta$,
but we do address more general cases in our experiments.

The articles Klein and Berrut \cite{BerrutKleinCAM} and Klein \cite{Klein} also claim that the
Lebesgue constant  of Klein's extended Floater-Hormann interpolant grows
logarithmically with $n$ and $d$, regardless of   $\tilde{n}$ and $\tilde{d}$.
Theorem 5.1 in page 6 of \cite{BerrutKleinCAM} summarizes this and other claims
from \cite{Klein} (see, in particular, Theorems 2.1 and 3.1 of \cite{Klein} and
the remark following the latter.) We cite:\\[0.02cm]

{\bf Theorem 5.1}(from \cite{BerrutKleinCAM})
{\em
Suppose $n$, $d$, $\tilde{n}$ and $\tilde{d}$ are positive integers, $\tilde{d} \leq \tilde{n} < n$
and assume that
$f \in C^{d+2}[a - d h, d + d h] \cap \wfc{C^{2 d + 1}}{[a, a + \tilde{n} h] \cup [b - \tilde{n} h,b]}$
is sampled at $n + 1$ equispaced nodes in $[a,b]$. Then
\begin{itemize}
\item[(i)] $\tilde{r}_n[f]$ has no real poles;
\item[(ii)]  For a constant $K$ independent of $n$, $
\wnorm{\tilde{r}_n[f] - f} \leq K h^{\min\wset{d,\tilde{d}} + 1}$;
\item[(iii)] The associated Lebesgue constant $\tilde{\Lambda}_n$ grows logarithmically with $n$ and $d$:
\[
\tilde{\Lambda}_n \leq 2 + \wfc{\ln}{n + 2 d}.
\]
\end{itemize}
}
Here we show that, in general, the third item in this theorem is false
if, as in page 2 of Berrut and Klein \cite{BerrutKleinCAM},  we assume the standard definition
of the Lebesgue constant as the norm of the interpolation operator,
and take into account the unavoidable errors in the extended function values $\tilde{y}$.
We prove that the traditional Lebesgue constant grows exponentially with $d$
when  $d = \tilde{n} = \tilde{d}$ and $\tilde{\wvec{y}}$ is outlined in \cite{BerrutKleinCAM}.
In the version of Theorem 5.1 stated in Klein \cite{Klein}
the reader is informed that actually the logarithmic bound assumes a peculiar interpretation of the Lebesgue
constant, namely, essentially that the mentioned approximate function values have no errors;
see the paragraph above Theorem 3.1 in \cite{Klein}. However, this
limitation of the result (iii) is not mentioned in \cite{BerrutKleinCAM}, and neither
\cite{BerrutKleinCAM} nor \cite{Klein} point out that the theorem does not apply
to the choices of $\tilde{\wvec{y}}$ proposed for the extended Floater-Hormann interpolants
and used in the experiments.

Rigorously, our proof applies only to the case $d = \tilde{d} = \tilde{n}$.
It suffices as a counterexample to Theorem 5.1, but it is unsatisfactory from
a broader practical perspective. However, we emphasize that,
in practice, Theorem 5.1 gives a misleading impression regarding
the Lebesgue constant of Extended Floater-Hormann interpolants
for broader classes of parameters. We do not have a formal theory
supporting this claim, but Sections \ref{section_practice} and
\ref{section_reduced} and Figure \ref{figure_conjecture} below
present strong experimental evidence of its validity.

\begin{figure}[!h]
\includegraphics[bb=-300 -10 260 340, width=6.0cm, height=3.2cm]{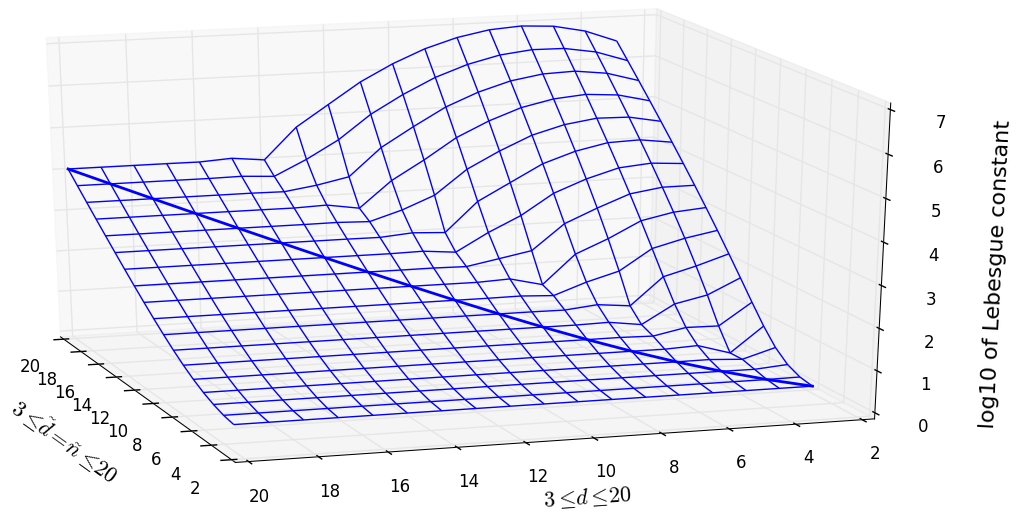}
\caption{Log10 of the correct Lebesgue constant as a function of $3 \leq d, \tilde{d} \leq 20$,
for $n = 100$ and $\tilde{n} = \tilde{d}$. Note that the diagonal $d = \tilde{d}$ highlighted in this figure
crosses the lines of constant $\tilde{d}$ in places in
which the Lebesgue constant is near the minimal value along such lines.
Therefore, it makes little sense to choose $d$ much less than $\tilde{d}$,
and the case considered in our counterexample is quite relevant in practice.
We have observed similar pictures for other values of $n$ and our experiments
indicate that the Lebesgue constant increases as $\tilde{n}$ gets larger than $\tilde{d}$.}
\label{figure_conjecture}
\end{figure}

This article presents an stability analysis of extended
Floater-Hormann interpolants based on the appropriate interpretation of the Lebesgue
constant. Formally, when the function $Y$ which yields the extended function values $\tilde{\wvec{y}}$ is
linear in $\wvec{y}$, we consider the linear operator $\wvec{y} \mapsto r_{\wvec{x},d,Y}[\wvec{y}] \in C^0[a,b]$
given by
\[
\wfc{r_{\wvec{x},d,Y}[\wvec{y}]}{t} := \wfc{\tilde{r}_{d, Y}}{t, \wvec{x}, \wvec{y}},
\]
for $\tilde{r}_{d, Y}$ defined in \wref{defrdy},
and the Lebesgue constant can be defined either as the norm of this linear
operator with respect to the supremum norm in $\wrn{n+1}$ and $C^0[x_0,x_n]$, or as the supremum of the Lebesgue function
\begin{equation}
\label{defleby}
\wfc{\Lambda_{\wvec{x}, d,Y}}{t} := \sup_{\wvec{y} \ \wrm{with} \ \wnorm{\wvec{y}}_\infty = 1 } \wabs{\wfc{r_{\wvec{x},d,Y}[\wvec{y}]}{t}}
\end{equation}
in $[x_0,x_n]$. These two definitions are equivalent, and
lead to a concept which has a fundamental
role  in theory and in practice, provided that it is
interpreted correctly.

By considering the correct Lebesgue constant, we gain
a more realistic view of the stability of extended Floater-Hormann
interpolants. For instance, we learn that the
version of these interpolants mentioned in Klein and Berrut
\cite{BerrutKleinCAM} and Klein \cite{Klein}
should not be used with large  $\tilde{d}$.
Since the order of approximation of these interpolants is
$h^\mu$ for $\mu =\min \wset{d,\tilde{d}}$, this limits their accuracy
in practice.

The articles \cite{BerrutKleinCAM} and \cite{Klein}
say nothing about the disastrous effect that a large $\tilde{d}$
may have on extended interpolants. Instead, they
emphasize that these interpolants can be used with large $d$.
This is illustrated by Figures 6 in  \cite{BerrutKleinCAM} and \cite{Klein},
which explore only the case $\tilde{d} = 7$ and compare the resulting extended
interpolants with usual interpolants with $d = \delta$ as large as $50$.
If instead of considering $d = \delta$ as large as 50 in their Figure 6,
they had focused on the more modest case $1 \leq d = \delta \leq \tilde{d} = 7$,
as in their other experiments, then they would have a more
realistic argument in favor of extended interpolants.
Indeed, for this range of $d = \delta$,
Figures 5 to 10 in \cite{Klein} show that extended interpolants are better
than usual ones in cases of practical interest. Therefore,
extended interpolants have merit and are a relevant topic for research.
However, our experiments
in Section \ref{section_instability} show that there
are also cases in which extended interpolants are worse than
usual ones, and here we aim at a balanced
view of their properties and limitations. In particular, we discuss
the role played by each one of their parameters and the ranges
in which they should be used.

In sections \ref{section_reduced} and \ref{section_lebesgue}
we analyze the Lebesgue constants of extended
interpolants from a theoretical perspective. We present an exponential
lower bound on the Lebesgue constant when $d = \tilde{n} = \tilde{d}$.
We also present experimental data showing that the dependency of the
Lebesgue function on these three parameters is not accurately
described by Theorem 5.1 in more general settings.
Section \ref{section_backward} discusses the backward stability of
extended Floater-Hormann interpolants in the general case in which
the function $\wfc{Y}{d,\wvec{x},\wvec{y}}$ that
defines $\tilde{\wvec{y}}$ is linear in $\wvec{y}$. We present
a simple condition that implies the backward instability of the barycentric formula
used to implement extended interpolants in this case, and we show experimentally that this
condition for backward instability is satisfied by an extended
interpolant mentioned in \cite{Klein}.

Section \ref{section_instability}
presents an empirical analysis of the stability of the extended
interpolants outlined in Berrut and Klein \cite{BerrutKleinCAM}.
We explain that the extrapolation
step may lead to numerical instability, and due to this instability the
overall error incurred by these interpolants can be much larger than
$\epsilon n \wleby{} \, \wnorm{\wvec{y}}_\infty$, where $\wleby{}$ is its Lebesgue constant and
$\epsilon$ is the machine precision.
In order to illustrate this fact, we present the results of
experiments in which the accuracy of the extended interpolants is
much worse than the accuracy of the usual interpolants.

On the positive side, once we become aware of the problems
caused by the extrapolation step, we may consider
ways to reduce them. When evaluating the interpolants for many
values of $t \in [x_0,x_n]$, it is worth computing the relatively few
extrapolated function values in multiple precision.
Numerical experiments show that this strategy leads to more accurate
extended interpolants.

Finally, the appendix considers the difficulties involved in the
construction of a general stability theory for  extended interpolants.
This appendix is at a more abstract level than the rest of the article:
we argue about the arguments one would use to discuss the stability of extended interpolants.
We hope that people interested in an in depth
analysis of the stability of these interpolants will appreciate
our remarks regarding the difficulties in formulating realistic hypotheses
and theorems about this subject.

\section{The rounding errors and the Lebesgue constant in practice}
\label{section_practice}
The Lebesgue constant is a fundamental concept in approximation theory.
It is also fundamental in practice, because it measures the sensitivity
of the interpolants to perturbations (or noise) in the data. In its proper interpretation,
the Lebesgue constant is equivalent to what numerical analysts call {\it condition number}, and
use to evaluate the numerical stability of algorithms.

\begin{figure}[!h]
\subfloat{
\includegraphics[bb=-60 10 320 187, width=4.7cm, height=2.5cm]{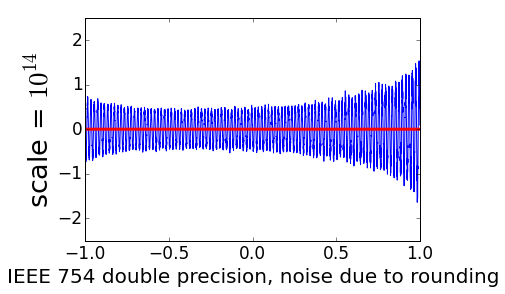}
}
\subfloat{
\includegraphics[bb=-130 20 620 260, width=7.7cm, height=2.5cm]{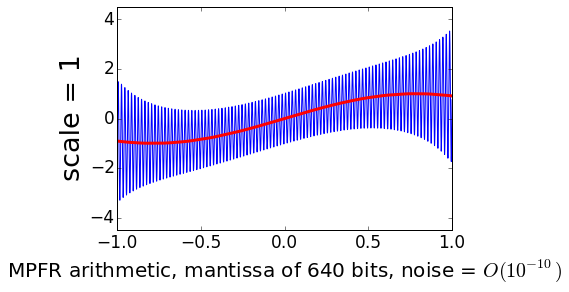}
}
\caption{The function $\wfc{f}{t} = \wfc{\sin}{2 t}$ (in red) and the approximation
obtained following the procedure proposed in \cite{BerrutKleinCAM}
with $n = 200$ (in blue). We use $d = \tilde{n} = \tilde{d} = 40$, whereas
\cite{BerrutKleinCAM} and \cite{Klein} consider $0 \leq d \leq 50$ and smaller values of
$\tilde{n}$ and $\tilde{d}$ in their experiments.
The  $\tilde{d}$ and $\tilde{n}$ in this figure satisfy the hypothesis of
Theorem 5.1 and $d$ is within the range considered in the experiments in \cite{BerrutKleinCAM} and \cite{Klein}.}
\label{figure_practice}
\end{figure}

This section shows that that rounding errors and noisy data have devastating effects on interpolants for
which theorems in \cite{BerrutKleinCAM} and \cite{Klein} claim that the Lebesgue
constant is small. Therefore, such claims may lead readers to believe that
these interpolants are much less affected by noise
and rounding errors than they really are.
Figure \ref{figure_practice} considers the approximation of $\wfc{f}{t} = \wfc{\sin}{2t}$
for $t \in [-1,1]$.  The plot on the left of Figure \ref{figure_practice} shows that by implementing the extrapolation
procedure proposed in \cite{BerrutKleinCAM} and \cite{Klein} with
the usual IEEE 754 double precision arithmetic we
may have numerical errors of order $10^{14}$ in circumstances
in which Theorem 5.1 yields a  bound smaller than $8$ on the constant
which Berrut and Klein call by Lebesgue's name.
These errors are several orders of magnitude larger than the ones reported in
\cite{BerrutKleinCAM} and \cite{Klein} for the same kind of extended
interpolant, because we do not restrict ourselves to the same small values
of $\tilde{d}$ as \cite{BerrutKleinCAM}
and \cite{Klein}.

The plot on the right of Figure \ref{figure_practice} illustrates the sensitivity of extended interpolants to
noise in the function values. It was obtained by adding random values of order
$10^{-10}$ to $y_0, \dots, y_n$. The effects of rounding errors in this plot
are negligible because we used the high precision arithmetic
provided by the MPFR library \cite{MPFR}, with a mantissa of 640 bits.
The experiment on the right indicates that in this case
the condition number is about $10^{10}$, and not $8$ as suggested by
Theorem 5.1 of Berrut and Klein \cite{BerrutKleinCAM}.

Figure \ref{figure_practice} raises
an interesting question: why does the plot on the
 left display errors of order
$10^{14}$ while the plot on the right shows errors of order one?
This question is intriguing because the IEEE 754's double precision machine
epsilon is of order $10^{-16}$ and is much smaller than the
$\wfc{O}{10^{-10}}$ perturbations used to generate the plot on the right.
As we explain in the rest of the article, the answer to this
question lies in the instabilities in the extrapolation process proposed
by Klein \cite{Klein} and this is one more reason why, in practice, the logarithmic bound
presented in \cite{BerrutKleinCAM} and \cite{Klein}
underestimates the effects of noise and rounding errors.

\section{The barycentric and reduced forms and the Lebesgue function}
\label{section_reduced}
In this section we show how to write extended interpolants in barycentric form
and introduce another way to describe them, which we call {\it reduced form}.
This form is numerically unstable and we do not advocate its use in practice.
Its purpose is to help us to deduce an expression for the Lebesgue function
of extended interpolants.

The Lebesgue function measures the sensitivity of the output of the complete
interpolation process to perturbations in its input, and we emphasize that the
input to the interpolation process are the original function
values $y_i$; not the extrapolated function values $\tilde{y}_i$.
Therefore, we can not ignore how changes in $\wvec{y}$ affect $\wvec{\tilde{y}}$
as suggested by Equation (3.2) in \cite{Klein}.

We recall that extended Floater-Hormann interpolants are defined only
for equally spaced nodes, and in this section focus on the interpolants with $\tilde{\wvec{y}}$ as
in the fifth section of \cite{BerrutKleinCAM}, ie.,  $\tilde{\wvec{y}}$ is defined using
extrapolation.
When the nodes are equally spaced,  \cite{Floater} shows that usual Floater-Hormann interpolants
can be written in the barycentric form
\begin{equation}
\label{bary_fh}
\wfh{t} = \ \ \left.
\sum\limits_{i = 0}^{n} \frac{w_{n,\delta,i} \, y_i}{t - x_i}
\right/ \sum\limits_{i = 0}^{n} \frac{w_{n,\delta,i}}{t - x_i},
\end{equation}
with weights
\begin{equation}
\label{bary_w}
w_{n,\delta,i} = (-1)^{i-\delta} \sum_{j = \max\{0, \, i - \delta\}}^{\min\{n-\delta, \, i\}}
\left( \begin{array}{c} \delta \\ i-j \end{array} \right),
\end{equation}
where the $y_i$ are the interpolated function values.
In the last paragraph of page 5 of \cite{BerrutKleinCAM},
extended interpolants are defined by extrapolating $\wvec{y}$ according to the following
Taylor series, which are defined in terms of the parameters $\tilde{d}$ and $\tilde{n}$:
\begin{eqnarray}
\label{ytila}
\tilde{y}_i & := & y_0 + \sum_{k=1}^{\tilde{d}} \wfc{r_{\underline{\wvec{x}},\tilde{d}}^{(k)}[\underline{\wvec{y}}]}{x_0}
\frac{\wlr{\tilde{x}_i - x_0}^k}{k!} \ \hspace{1.5cm} \wrm{for} \ -d \leq i < 0, \\
\label{ytilb}
\tilde{y}_i & := & y_i \hspace{5.82cm} \wrm{for} \ 0 \leq i \leq n, \\
\label{ytilc}
\tilde{y}_i & := & y_n + \sum_{k=1}^{\tilde{d}} \wfc{r_{\overline{\wvec{x}},\tilde{d}}^{(k)}[\overline{\wvec{y}}]}{x_n}
\frac{\wlr{\tilde{x}_i - x_n}^k}{k!} \ \hspace{1.4cm} \wrm{for} \ n < i \leq n + d,
\end{eqnarray}
where $\tilde{x}_i = x_0 + i h$ for $-d \leq i \leq n + d$  and $\wfhk{t}$ is the $k$th derivative
of the Floater-Hormann interpolant $\wfc{r_{\wvec{x},d}[\wvec{y}]}{t}$ in \wref{bary_fh}
and $\underline{\wvec{x}} := \wlr{x_0,\dots, x_{\tilde{n}} }$,
$\underline{\wvec{y}} := \wlr{y_0,\dots, y_{\tilde{n}} }$,
$\overline{\wvec{x}} := \wlr{x_{n - \tilde{n}},\dots, x_n}$ and
$\overline{\wvec{y}} := \wlr{y_{n - \tilde{n}},\dots, y_n}$.
Therefore, the extended interpolant is specified once we define $\wvec{x}$, $d$, $\tilde{n}$ and $\tilde{d}$,
and we can write it in the following barycentric form:
\begin{equation}
\label{bary_xfh}
\wxfh{t} = \ \ \left.
\sum\limits_{i = -d}^{n + d} \frac{\tilde{w}_{n,d,i} \, \tilde{y}_i}{t - \tilde{x}_i} \right/ \sum\limits_{i = -d}^{n + d} \frac{\tilde{w}_{n,d,i}}{t - \tilde{x}_i},
\end{equation}
with weights
\[
\tilde{w}_{n,d,i} := w_{n + 2 d,d,i + d} = (-1)^{i} \sum_{j = \max\{0, \, i\}}^{\min\{n  +d, \, i + d\}} \left( \begin{array}{c} d \\ i + d -j \end{array} \right).
\]
It is difficult to derive
the Lebesgue function of the extended interpolant $\wxfhi{}$ directly from Equation \wref{bary_xfh},
because this equation depends on $\tilde{\wvec{y}}$, which is not part of the
original interpolation problem. To derive an expression for this Lebesgue function, it
is helpful to write the extended interpolant only in terms of the original $\wvec{y}$.
The next lemma explains how to achieve this goal when $2 \tilde{n} < n$:
\begin{lemma}
\label{lem_clean}
Given extended nodes $\tilde{x}_i = x_0 + i h$ for $-d \leq i \leq n + d$,
the extrapolated function values $\tilde{y}_i$ used to define the extended interpolant $\wxfhi{}$ can be written as
\begin{eqnarray}
\label{tilyi_clean_a}
\tilde{y}_i & = & \sum_{j = 0}^{\tilde{n}} a_{ij} y_j,   \hspace{2.25cm} \wrm{for} \hspace{0.1cm} - d \leq i < 0, \\
\label{tilyi_clean_b}
\tilde{y}_i & = & \sum_{j = n - \tilde{n}}^n b_{\wlr{i-n}\wlr{j - n}} y_j, \hspace{0.6cm} \wrm{for} \hspace{0.1cm} n < i \leq n + d,
\end{eqnarray}
where the numbers
\[
\wset{a_{ij}, \ -d \leq i < 0, \ 0 \leq j \leq \tilde{n}} \hspace{0.5cm} \wrm{and} \hspace{0.5cm}
\wset{b_{ij}, \ 0 < i \leq d, \ - \tilde{n} \leq j \leq 0}
\]
depend on $i$,$j$,$\tilde{n}$ and $\tilde{d}$ but do not depend on $h$, $d$ or $n$, in the sense that
there exist functions $\alpha, \beta: \wz{}^4 \rightarrow \wrone{}$ such that $a_{ij} = \wfc{\alpha}{i,j,\tilde{n},\tilde{d}}$
and $b_{ij} = \wfc{\beta}{i,j,\tilde{n},\tilde{d}}$.
When $2 \tilde{n} < n$ the extended interpolant can be written in the reduced form
\begin{equation}
\label{reduced_form}
\wxfh{t} = \left. \sum_{j = 0}^n \wfc{c_j}{t} y_j \right/ \sum_{j = -d}^{n+d} \frac{\tilde{w}_{n,d,j}}{t - \tilde{x}_j},
\end{equation}
where the functions $c_j$ are given by
\begin{eqnarray}
\label{ca}
\hspace{0.4cm}
\wfc{c_{j}}{t} & = & \frac{\tilde{w}_{n,d,j}}{t - x_{j}} + \sum_{i = -d}^{-1} \frac{\tilde{w}_{n,d,i} a_{ij}}{t - \tilde{x}_{i}}
\hspace{2.35cm} \wrm{for} \ 0 \leq j \leq \tilde{n},\\
\label{cb}
\hspace{0.4cm}
\wfc{c_{j}}{t} & = & \frac{\tilde{w}_{n,d,j}}{t - x_{j}}
\hspace{4.9cm} \wrm{for} \ \tilde{n} < j < n - \tilde{n}, \\
\label{cc}
\hspace{0.4cm}
\wfc{c_{j}}{t} & = & \frac{\tilde{w}_{n,d,j}}{t - x_{j}} + \sum_{i = n + 1}^{n + d} \frac{\tilde{w}_{n,d,i} b_{\wlr{i-n}\wlr{j-n}}}{t - \tilde{x}_{i}}
\hspace{1.0cm} \wrm{for} \ n - \tilde{n} \leq j \leq n .
\end{eqnarray}
\end{lemma}
In the end of this section we prove Lemma \ref{lem_clean} and present explicit expressions for
$a_{ij}$ and $b_{ij}$. We also provide formulae analogous to
\wref{ca}--\wref{cc} for the case $2 \tilde{n} \geq n$.

The Lebesgue functions of the interpolants $\wfhi{}$ and $\wxfhi{}$ are defined as
\[
\wlebffh{t} := \sup_{\wnorm{\wvec{y}}_\infty = 1} \wabs{ \wfh{t}}
\hspace{1cm} \wrm{and} \hspace{1cm}
\wlebfxfh{t} := \sup_{\wnorm{\wvec{y}}_\infty = 1} \wabs{ \wxfh{t}},
\]
and using Equation \wref{bary_fh} and the reduced form \wref{reduced_form} it is easy to show that
\begin{equation}
\label{leb_fun_fh}
\wlebffh{t} =
\ \ \left.
\sum\limits_{j = 0}^{n} \wabs{\frac{w_{n,\delta,j}}{t - x_j}} \right/ \wabs{\sum\limits_{j = 0}^{n} \frac{w_{n,\delta,j}}{t - x_j}}
\end{equation}
and
\begin{equation}
\label{leb_fun_xfh}
\wlebfxfh{t} =  \ \
\left. \sum_{j = 0}^n \wabs{\wfc{c_j}{t} } \right/ \wabs{\sum_{j = -d}^{n+d} \frac{\tilde{w}_{n,d,j}}{t - \tilde{x}_j}}.
\end{equation}
We emphasize that, in general,
\begin{equation}
\label{wrong}
\wlebfxfh{t} \neq
\ \left. \sum\limits_{j = -d}^{n + d} \wabs{\frac{\tilde{w}_{n,d,j}}{t - \tilde{x}_j}}
\right/ \wabs{\sum\limits_{j = -d}^{n + d} \frac{\tilde{w}_{n,d,j}}{t - \tilde{x}_j}},
\end{equation}
that is, we can deduce \wref{leb_fun_fh} from \wref{bary_fh}, but
$\wlebfxfh{t}$ is not equal to, or even bounded by, the right hand side of
\wref{wrong}. This is why Equation \kleinsbug{} in \cite{Klein} is misleading.
The fact that this equation refers to the peculiar interpretation
of the Lebesgue constant used in \cite{Klein}, and not to the actual Lebesgue constant, becomes
evident when we plot the right and the left hand side of \wref{wrong}  for
 $n = 50$, $d = 3$, $\tilde{n} = 11$ and $\tilde{d} = 7$, as in Figure \ref{figure_kleins_blunder} below.

\begin{figure}[!h]
\subfloat[The right hand side of Equation \kleinsbug{} in \cite{Klein} as a function of $t \in {[}-1,1{]}$, which is
claimed to be an upper bound on the Lebesgue function]{
\includegraphics[bb=-60 10 320 240, width=5cm, height=2.1cm]{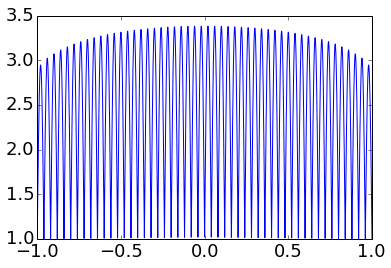}
}
\subfloat[The correct Lebesgue function]{
\includegraphics[bb=-130 10 520 240, width=8cm, height=2.1cm]{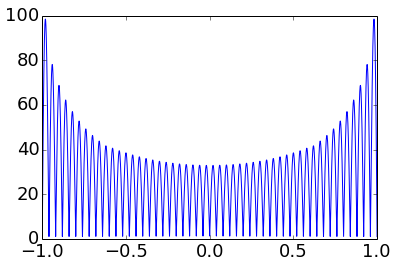}
}
\caption{The correct Lebesgue function for
$n = 50$, $d = 3$, $\tilde{n} = 11$ and $\tilde{d} = 7$.}
\label{figure_kleins_blunder}
\end{figure}

\begin{figure}[!h]
\includegraphics[bb=-80 0 1020 900, width=10cm, height=5.5cm]{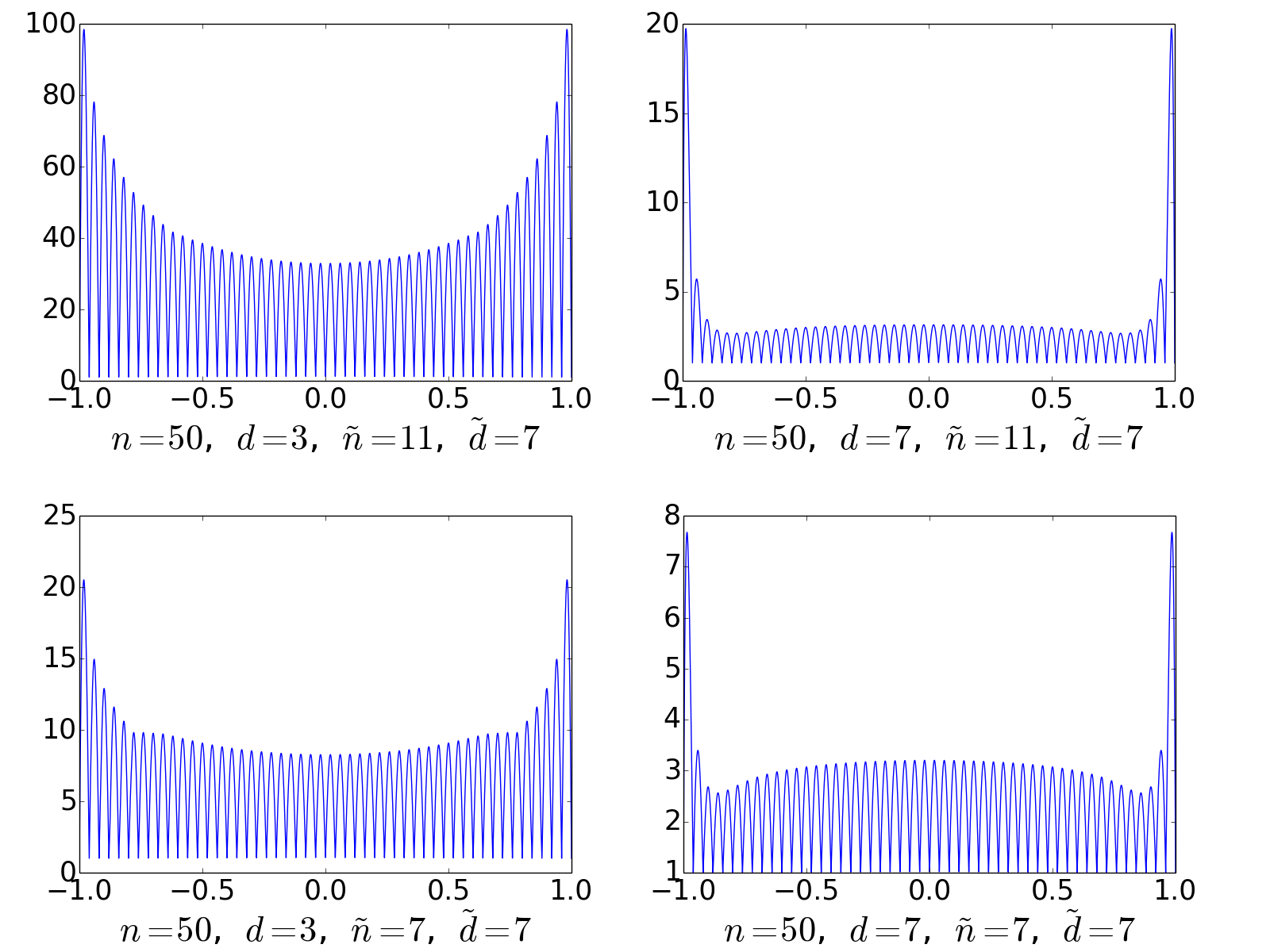}
\caption{The dependency of the Lebesgue function on $d$, $\tilde{n}$ and $\tilde{d}$.}
\label{figure_lebesgue_functions}
\end{figure}

The dependency of the Lebesgue function on the parameters $d$, $\tilde{n}$ and $\tilde{d}$
is subtle. For instance, Figure \ref{figure_lebesgue_functions} shows that
the Lebesgue function may decrease as we increase $d$ and keep the other parameters fixed.
Moreover, the Lebesgue constant for $n = 50$ and $d = \tilde{n} = \tilde{d} = 7$
is $12$ times smaller than the Lebesgue constant if
$n = 50,$  $d = 3$, $\tilde{n} = 11$ and $\tilde{d} = 7$, as
considered in \cite{Klein}.

\subsection{Proof of Lemma \ref{lem_clean}}
\label{subsection_lem_clean}
In this subsection we prove Lemma \ref{lem_clean} and show that when
$n \leq 2 \tilde{n}$ we have an analogous result with
\begin{eqnarray}
\label{cjba}
\hspace{0.4cm}
\wfc{c_{j}}{t} & = & \frac{\tilde{w}_{n,d,j}}{t - x_{j}} + \sum_{i = -d}^{-1} \frac{\tilde{w}_{n,d,i} \, a_{i j}}{t - \tilde{x}_{i}}
\hspace{2.7cm} \wrm{for} \ 0 \leq j < n - \tilde{n},   \\
\hspace{0.4cm}
\nonumber
\wfc{c_{j}}{t} & = & \frac{\tilde{w}_{n,d,j}}{t - x_{j}}  + \sum_{i = -d}^{-1} \frac{\tilde{w}_{n,d,i}\,  a_{i j}}{t - \tilde{x}_{i}}\\
\label{cjbb}
& + &  \sum_{i = n + 1}^{n + d} \frac{\tilde{w}_{n, d, i} \, b_{\wlr{i-n}\wlr{j-n}}}{t - \tilde{x}_i}
 \hspace{2.6cm} \wrm{for} \ n - \tilde{n} \leq j \leq \tilde{n}, \\
 \label{cjbc}
\wfc{c_{j}}{t} & = & \frac{\tilde{w}_{n,d,j}}{t - x_{j}} +
\sum_{i = n + 1}^{n + d} \frac{\tilde{w}_{n,d,i} \,  b_{\wlr{i-n} \wlr{j-n}}}{t - \tilde{x}_i}
\hspace{1.27cm} \wrm{for} \ \tilde{n} < j \leq n .
\end{eqnarray}
Our proof involves the numbers $D^{(k)}_{ij}$ mentioned in the third section of \cite{KleinBerrutDerivatives},
which represent the $k$th derivative at the node $x_i$ of the $j$th Lagrange fundamental rational function.
In order to simplify the notation, we use the following lemma.
\begin{lemma}
\label{lem_e}
Consider nodes $x_0 < x_1 < \dots < x_{\tilde{n}}$, weights $w_0, w_1,\dots, w_{\tilde{n}} \in \wrone{} - \wset{0}$, and the numbers
$E^{(k)}_{ij}$, with $k \geq 0$ and $0 \leq i, j \leq \tilde{n}$, defined inductively in $k$ by
\begin{eqnarray}
\label{eijkA}
E^{(0)}_{ii} & := & 1 \hspace{1cm} \wrm{and} \hspace{1cm} E^{(0)}_{ij} := 0 \hspace{1.42cm} \wrm{for} \ i \neq j. \\
\label{eijkB}
E^{(k)}_{ij} & := & \frac{k}{x_i - x_j} \wlr{\frac{w_j}{w_i} E_{ii}^{\wlr{k-1}} - E_{ij}^{\wlr{k-1}}}
\hspace{0.9cm} \wrm{for} \ i \neq j, \\
\label{eijkC}
E_{ii}^{(k)} & := & - \sum_{j = 0, j \neq i}^{\tilde{n}} E_{ij}^{(k)}.
\end{eqnarray}
If $k > 0$ then $E^{(k)}_{ij}$ is equal to the number $D^{(k)}_{ij}$ in Equations (3.1) and (3.2) in
\cite{KleinBerrutDerivatives} with $n$ replaced by $\tilde{n}$.
\end{lemma}
In order to prove this lemma, replace $k$ by $1$ in Equations \wref{eijkB} and \wref{eijkC},
use \wref{eijkA} to simplify the result and
Equation (3.2) of \cite{KleinBerrutDerivatives} to verify that
$E^{(1)}_{ij} = D^{(1)}_{ij}$. The case $k > 1$ follows
by induction from \wref{eijkB} and \wref{eijkC}  and
Equation (3.3) in \cite{KleinBerrutDerivatives}.

Since we are assuming that the nodes are equally spaced, $x_i - x_j = \wlr{i - j} h$ and it is convenient to consider
the normalized numbers
\[
\overline{E}_{ij}^{(k)} := h^k E_{ij}^{(k)} / k!.
\]
Replacing $E_{ij}^{(k)}$ by $\overline{E}_{ij}^{(k)} k! \, h^{-k}$ in \wref{eijkA}--\wref{eijkC}
we obtain
\begin{eqnarray}
\label{final_d0}
\hspace{0.4cm}
\overline{E}^{(0)}_{ii} & = & 1 \hspace{1cm} \wrm{and} \hspace{1cm} \overline{E}^{(0)}_{ii} = 0
\hspace{1.73cm} \wrm{for} \ i \neq j, \\
\label{final_d1}
\hspace{0.4cm}
\overline{E}^{(k)}_{ij} & = & \frac{1}{i- j} \wlr{\frac{w_j}{w_i} \overline{E}_{ii}^{\wlr{k-1}} - \overline{E}_{ij}^{\wlr{k-1}}}
\hspace{1.47cm} \wrm{for} \ i \neq j \hspace{0.2cm} \wrm{and} \hspace{0.2cm} k \geq 1, \\
\label{final_d2}
\hspace{0.4cm}
\overline{E}_{ii}^{(k)} & = & - \sum_{j = 0, j \neq i}^{\tilde{n}} \overline{E}_{i j}^{(k)}.
\end{eqnarray}
Equations \wref{final_d0}--\wref{final_d2} show that $\overline{E}^{(k)}_{ij}$ depends on $i$, $j$, $k$, and $\tilde{n}$,
but it can depend on $h$, $d$, $n$ or $\tilde{d}$ only via the weights $w_i$, because
there is no mention to $h$, $d$, $n$ and $\tilde{d}$ in \wref{final_d0} -- \wref{final_d2}.
In the case that concerns us, the weights $w_i$ correspond to the
usual Floater-Hormann interpolants in equally spaced nodes $x_0,\dots, x_{\tilde{n}}$ with
parameter $\delta = \tilde{d}$. These weights are given by \wref{bary_w},
with $n$  and $\delta$ replaced by $\tilde{n}$ and $\tilde{d}$,
and depend on $\tilde{n}$ and $\tilde{d}$, but not on $h$, $n$ or $d$.
Therefore, in the case relevant to our discussion, $\overline{E}^{(k)}_{ij}$ does not depend on
$h$, $n$ or $d$.

Equation (3.1) in \cite{KleinBerrutDerivatives} and the identities $D_{ij}^{(k)} = E^{(k)}_{ij}$ for $k > 0$
imply that
\begin{eqnarray}
\nonumber
\wfc{r^{(k)}_{\underline{\wvec{x}},\tilde{d}}[\underline{\wvec{y}}]}{x_0} & = & \sum_{j = 0}^{\tilde{n}} D_{0j}^{(k)} y_j =
\sum_{j = 0}^{\tilde{n}} E_{0j}^{(k)} y_j =
\frac{k!}{h^k} \sum_{j = 0}^{\tilde{n}} \overline{E}_{0j}^{(k)} y_j,  \\
\nonumber
\wfc{r^{(k)}_{\overline{\wvec{x}},\tilde{d}}[\overline{\wvec{y}}]}{x_n} & = &
\sum_{j = n - \tilde{n}}^{n} D_{\tilde{n}\wlr{j -n + \tilde{n}}}^{(k)} y_j =
\sum_{j = n - \tilde{n}}^{n} E_{\tilde{n}\wlr{j -n + \tilde{n}}}^{(k)} y_j
= \frac{k!}{h^k} \sum_{j = n - \tilde{n}}^{n} \overline{E}_{\tilde{n}\wlr{j - n + \tilde{n}}}^{(k)} y_j,
\end{eqnarray}
for $k > 0$. Combining the last two equations with the identities $\overline{E}^{(0)}_{ii} = 1$
and $\overline{E}^{(0)}_{ij} = 0$ for $i \neq j$
we can rewrite \wref{ytila} and \wref{ytilc} as
\begin{eqnarray}
\nonumber
\tilde{y}_i & = & \sum_{k=0}^{\tilde{d}} \sum_{j = 0}^{\tilde{n}} \overline{E}_{0j}^{(k)} i^k y_j
\ \hspace{3.4cm} \wrm{for} \ -d \leq i < 0, \\
\nonumber
\tilde{y}_i & = & \sum_{k=0}^{\tilde{d}} \sum_{j = n - \tilde{n}}^{n} \overline{E}_{\tilde{n}\wlr{j - n + \tilde{n}}}^{(k)} \wlr{i-n}^k y_j \
\hspace{1.0cm} \wrm{for} \ n < i \leq n + d.
\end{eqnarray}
These equations are equivalent to \wref{tilyi_clean_a} and \wref{tilyi_clean_b} with
\begin{eqnarray}
\label{def_aij}
\hspace{0.8cm} a_{ij}  & := & \sum_{k = 0}^{\tilde{d}}  \overline{E}_{0j}^{(k)} i^k
\hspace{1.3cm} \wrm{for}
\ -d \leq i < 0 \hspace{0.3cm} \wrm{and} \hspace{0.3cm} 0 \leq j \leq \tilde{n}, \\
\label{def_bij}
\hspace{0.8cm}  b_{i j } & := & \sum_{k = 0}^{\tilde{d}} \overline{E}_{\tilde{n} \wlr{j + \tilde{n}}}^{(k)} i^k
\hspace{0.8cm} \wrm{for}
\ 0 < i \leq d \hspace{0.3cm} \wrm{and} \hspace{0.3cm} - \tilde{n} \leq j \leq 0.
\end{eqnarray}
Since $\overline{E}_{ij}^{(k)}$ does not depend on $h$, $d$ or $n$,
and there is no mention to $h$, $d$ or $n$ in the right hand side of
Equations \wref{def_aij} and \wref{def_bij},
its is clear that $a_{ij}$ and $b_{ij}$ do not depend on $h$, $d$ or $n$,
as claimed in Lemma \ref{lem_clean}.
Equation \wref{bary_xfh} yields
\[
\wxfh{t} = \frac{1}{\wfc{Q}{t}} \sum_{i = -d}^{n + d} \frac{\tilde{w}_{n,d,i} \, \tilde{y}_i}{t - \tilde{x}_i}
\hspace{1.0cm} \wrm{with} \hspace{1.0cm}
\wfc{Q}{t} = \sum_{i = -d}^{n + d} \frac{\tilde{w}_{n,d,i}}{t - \tilde{x}_i}.
\]
It follows that
\[
\wfc{Q}{t} \wxfh{t}
= \sum_{i = -d}^{-1} \frac{\tilde{w}_{n,d,i}}{t - \tilde{x}_i} \sum_{j = 0}^{\tilde{n}} \, a_{i j} \,  y_j  +
  \sum_{j = 0}^{n} \frac{\tilde{w}_{n,d,j} \, y_j}{t - x_j} \ +
\]
\[
  \sum_{i = n + 1}^{n + d} \frac{\tilde{w}_{n,d,i}}{t - \tilde{x}_i} \sum_{j = n - \tilde{n}}^n \,  b_{\wlr{i-n} \wlr{j-n}} \, y_j  =
\]
\begin{equation}
\label{qform}
 \sum_{j = 0}^{\tilde{n}} \sum_{i = -d}^{-1} \frac{\tilde{w}_{n,d,i} \,  a_{i j} \,  y_j}{t - \tilde{x}_i} +
 \sum_{j = 0}^{n} \frac{\tilde{w}_{n,d,j} \, y_j}{t - x_j} +
 \sum_{j = n - \tilde{n}}^n \sum_{i = n + 1}^{n + d} \frac{\tilde{w}_{n,d,i} \,  b_{\wlr{i-n}\wlr{j-n}} \,  y_j}{t - \tilde{x}_i}.
\end{equation}
We now have two cases: (i) $\tilde{n} < n - \tilde{n}$ and (ii) $\tilde{n} \geq n - \tilde{n}$.
In the first case we can rewrite \wref{qform} as
\[
\wfc{Q}{t} \wxfh{t} =
 \sum_{j = 0}^{\tilde{n}} \wlr{\frac{\tilde{w}_{n,d,j}}{t - x_j} +
 \sum_{i = -d}^{-1} \frac{\tilde{w}_{n,d,i} \, a_{i j}}{t - \tilde{x}_i} } y_j \ +
 \]
 \[
  \sum_{j = \tilde{n} + 1}^{n - \tilde{n} - 1} \frac{\tilde{w}_{n,d,j}}{t - x_j} \,  y_j +
  \sum_{j = n - \tilde{n}}^n \wlr{\frac{\tilde{w}_{n,d,j}}{t - x_j} +
  \sum_{i = n + 1}^{n + d} \frac{\tilde{w}_{n,d,i} \, b_{\wlr{i-n}\wlr{j-n}}}{t - \tilde{x}_i}} y_j,
\]
and this proves \wref{reduced_form}--\wref{cc}. When $\tilde{n} \geq n - \tilde{n}$,
we can rewrite \wref{qform} as
\[
\wfc{Q}{t} \wxfh{t} =
 \sum_{j = 0}^{n - \tilde{n} - 1} \wlr{\frac{\tilde{w}_{n,d,j}}{t - x_j} +
 \sum_{i = -d}^{-1} \frac{\tilde{w}_{n,d,i} \, a_{i j}}{t - \tilde{x}_i} } y_j \ +
 \]
 \[
  \sum_{j = n - \tilde{n}}^{\tilde{n}} \wlr{\sum_{i = -d}^{-1} \frac{\tilde{w}_{n,d,i} \,  a_{i j}}{t - \tilde{x}_i} +
  \frac{\tilde{w}_{n,d,j}}{t - x_j} +
\sum_{i = n + 1}^{n + d} \frac{\tilde{w}_{n,d,i} \, b_{\wlr{i-n}\wlr{j-n}}}{t - \tilde{x}_i}  } y_j \ +
 \]
 \[
 \sum_{j = \tilde{n} + 1}^n \wlr{\frac{\tilde{w}_{n,d,j}}{t - x_j} +
 \sum_{i = n + 1}^{n + d} \frac{\tilde{w}_{n,d,i} \, b_{\wlr{i-n}\wlr{j-n}}}{t - \tilde{x}_i}} y_j.
\]
This proves \wref{cjba}--\wref{cjbc} and we are done. \qed{}

\section{The Lebesgue constant when $d = \tilde{n} = \tilde{d}$}
\label{section_lebesgue}
This section presents a proof that the Lebesgue constant
of  extended interpolants mentioned in \cite{BerrutKleinCAM} and \cite{Klein} grows exponentially with $d = \tilde{n} = \tilde{d}$.
This shows that the peculiar interpretation of the
Lebesgue constant mentioned in \cite{Klein}
does not capture essential points regarding the stability of
extended Floater Hormann interpolants in general, because
Equation (3.2) in \cite{Klein}  does not
take properly into account how changes on $\wvec{y}$ affect $\tilde{\wvec{y}}$.

The Lebesgue constant of the extended interpolant $\wxfh{t}$ in \wref{bary_xfh}
is
\[
\wlebcxfh{} := \max_{x_0 \leq t \leq x_n} \wlebfxfh{t},
\]
for the Lebesgue function $\wlebfxfh{t}$ in \wref{leb_fun_xfh}, and the Lebesgue constant for polynomial interpolation at $d + 1$ equally
spaced nodes is
\[
\Lambda_d := \sup_{0 \leq t \leq d, \ \ \wvec{y} \in \wrn{d+1} - \wset{0}} \frac{\wabs{ \wfc{p_{d,\wvec{y}}}{t}}}{\wnorm{\wvec{y}}_\infty},
\]
where $p_{d,\wvec{y}}$ is the polynomial with degree less than $d + 1$ such that $\wfc{p_{d,\wvec{y}}}{i} = y_i$ for $i = 0,\dots,d$.
In this section we show that $\wlebcdiag{}$ is not much smaller than $\Lambda_d$,
by providing a lower bound for $\wlebcdiag{}$ which
approaches $\Lambda_d$ exponentially fast as $d$ increases.
Formally, we have the following:

\begin{theorem}
\label{thm_main}
If $n > d + 1 \geq 3$ then $\wlebcdiag{} \geq   \kappa_d  \Lambda_d$, for
\begin{equation}
\label{main_thesis}
\kappa_d :=  \wlr{1 - \frac{d}{2^d - 1} - 2^{-d}}.
\end{equation}
\end{theorem}

\begin{figure}[!h]
\includegraphics[bb=0 0 1100 650, width=10cm, height=4.6cm]{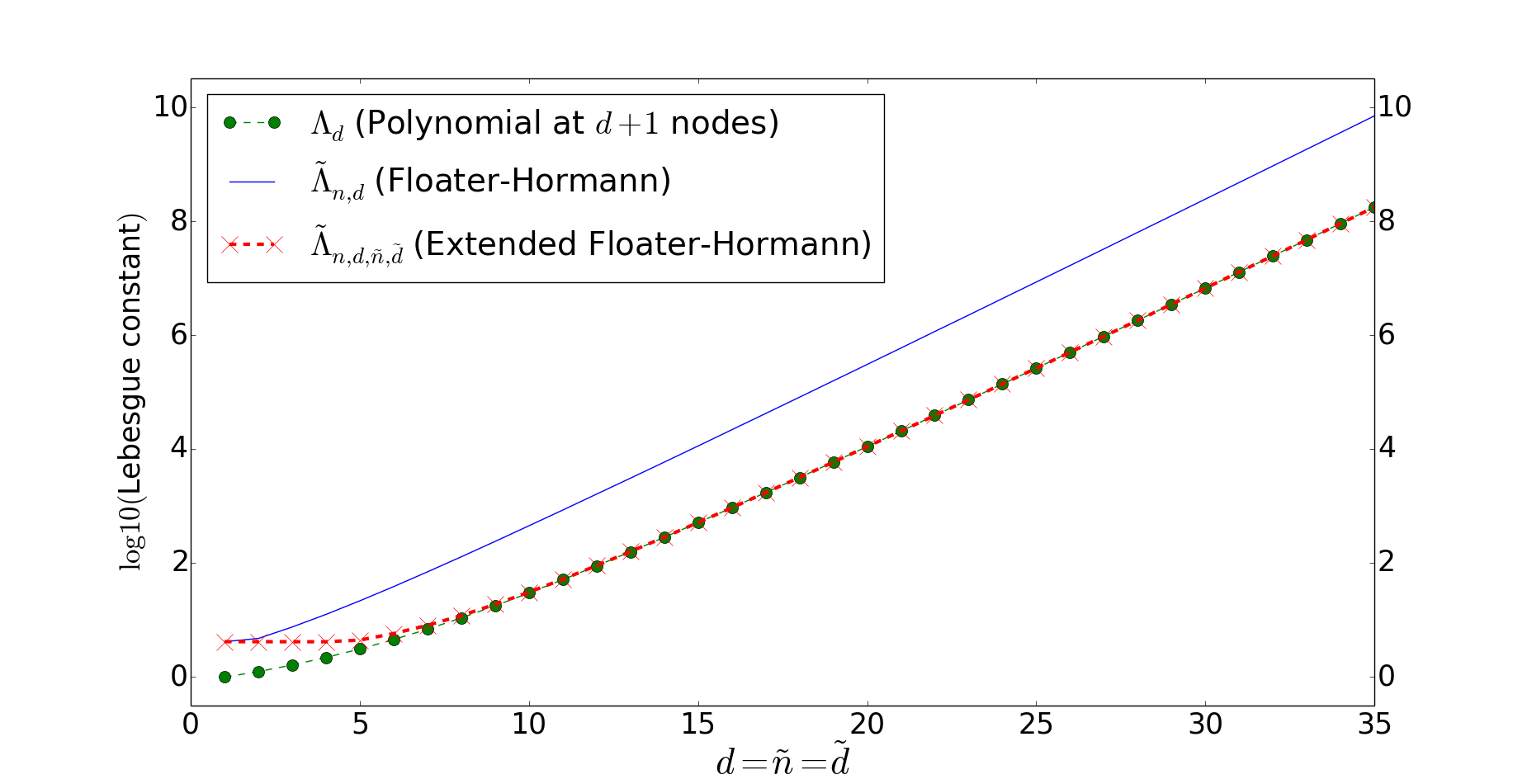}
\caption{
$\log\!{\rm10}$ of the Lebesgue constants for $n = 200$ and
$\delta = d = \tilde{n} = \tilde{d}$ varying from $1$ to $35$. The Lebesgue constant of the extended interpolant
is about the same as the Lebesgue constant for polynomial interpolation at $\wlr{d+1}$ equally
spaced nodes for $d \geq 7$, and $\wlebcfh{}$ is roughly equal to $100 \wlebcxfh{}$
for large $\delta = d = \tilde{n} = \tilde{d}$.
}
\label{figure_lebesgue_constants}
\end{figure}

We prove Theorem \ref{thm_main} at the end of this section. For now, let us
explore its consequences and check them experimentally.
As explained in \cite{Trefe}, we have
\[
\frac{2^{d-2}}{d^2} < \Lambda_d < \frac{2^{d+3}}{d}.
\]
Therefore, $\Lambda_d$ grows exponentially with $d$ and Theorem \ref{thm_main} shows that
the same applies to $\wlebcdiag{}$. Moreover, Bos et. al. \cite{Bos} show that
the Lebesgue constant for the Floater-Hormann interpolant
at $\wvec{x}$ with parameter $\delta$ satisfies
$\wlebcfh{} \leq 2^{\delta-1} \wlr{2 + \log n}$.
Theorem \ref{thm_main} shows that $\Lambda_d < 1.5 \wlebcdiag{}$ for $d \geq 4$,
and combining the two equations above for $d \geq 4$ we conclude that,
when $\delta = d$,
\[
\wlebcfh{} < 2 d^2 \wlr{2 + \log n} \Lambda_d < 3 d^2 \wlr{2 + \log n} \wlebcdiag{},
\]
and the ratio $\wlebcfh{}/\wlebcdiag{}$ is definitely not
as large as claimed in \cite{Klein}. This observation is
corroborated by Figure \ref{figure_lebesgue_constants}.

%%%%%%%%%%%%%%%%%%%%%%%%%%%%%%%%%%%%%%%%%%%%%%%%%%%%%%%%%%%%%%%%%%%%%%%%%%%%%%%%%%%%%%%%%%%%
%%%%%%%%%%%%%%%%%%%%%%%%%%%%%%%%%%%%%%%%%%%%%%%%%%%%%%%%%%%%%%%%%%%%%%%%%%%%%%%%%%%%%%%%%%%%

\subsection{Proof of Theorem \ref{thm_main}}
We follow the usual convention that a sum of the form $\sum_{i = a}^b u_i$ with $b < a$ is equal
to $0$ and a product $\prod_{i = a}^b u_i$ with $b < a$ is equal to $1$. According to \cite{Klein},
\begin{equation}
\label{FHinterp2}
\wxfh{t}  = \frac{\sum\limits_{i = -d}^{n} \tilde{\lambda}_i(t,\tilde{\wvec{x}})
  p_i(t,\tilde{\wvec{x}},\wvec{\tilde{y}}) }{\sum\limits_{i = -d}^{n} \tilde{\lambda}_i(t,\tilde{\wvec{x}})  },
\end{equation}
where $\tilde{x}_i = x_0 + i h$ for $-d \leq i \leq n + d$,
$\tilde{\wvec{y}} = \wlr{\tilde{y}_{-d},\tilde{y}_{1-d}, \, \dots \,, \tilde{y}_{n+d}}$
and
\begin{equation}
\label{LambdaAnd}
\wfc{\tilde{\lambda}_i}{t,\tilde{\wvec{x}}} := \frac{\wlr{-1}^i}{\wlr{t - \tilde{x}_i} \wlr{t - \tilde{x}_{i+1}}
\dots \wlr{t - \tilde{x}_{i+d}}} \hspace{1cm} \wrm{for} \ i = 0, \dots, n,
\end{equation}
and $\wfc{p_{i}}{t,\tilde{\wvec{x}},\tilde{\wvec{y}}}$ is the polynomial with degree less than $d + 1$ such that
$\wfc{p_{i}}{\tilde{x}_k,\tilde{\wvec{x}},\tilde{\wvec{y}}} = \tilde{y}_k$ for $k = i,\dots,i + d$.
When $\tilde{\wvec{y}}$ is defined as in Equations \wref{ytila}-- \wref{ytilc}, we have
\[
\wfc{p_i}{t,\tilde{\wvec{x}},\wvec{\tilde{y}}} = \wfc{p_i}{t,\wvec{x},\wvec{y}}  \ \ \ \wrm{for} \ \ \ 0 \leq i \leq n - d,
\]
and when $d = \tilde{n} = \tilde{d}$ we also have
\begin{eqnarray}
\nonumber
\wfc{p_i}{t,\tilde{\wvec{x}},\wvec{\tilde{y}}} & = & \wfc{p_0}{t,\wvec{x},\wvec{y}}  \hspace{1.4cm} \ \wrm{for} \  -d \leq i < 0,\\
\nonumber
\wfc{p_i}{t,\tilde{\wvec{x}},\wvec{\tilde{y}}} & = &  \wfc{p_{n-d}}{t,\wvec{x},\wvec{y}} \hspace{1cm} \ \wrm{for} \   n-d < i \leq n,
\end{eqnarray}
because in this case the interpolants $\wfc{r_{\underline{\wvec{x}},\tilde{d}}[\underline{\wvec{y}}]}{t}$ and
$\wfc{r_{\overline{\wvec{x}},\tilde{d}}[\overline{\wvec{y}}]\dot{}}{t}$ are polynomials,
and the Taylor series of a polynomial is equal to itself.
Equation \wref{FHinterp2} leads to
\[
\wxfhdiag{t} =
 \frac{
 \left(\sum_{i = -d}^{0} \wfc{\tilde{\lambda}_i}{t,\tilde{\wvec{x}}}\right)
 \wfc{p_0}{t,\wvec{x},\wvec{y}}
 +\wfc{\tilde{\lambda}_1}{t,\wvec{x}} \wfc{p_1}{t,\wvec{x},\wvec{y}}}
{\sum_{i = -d}^{n} \wfc{\tilde{\lambda}_i}{t,\tilde{\wvec{x}}}} \ +
\]
\begin{equation}
\label{FHinterp3}
  \frac{\sum_{i = 2}^{n-d - 1} \wfc{\tilde{\lambda}_i}{t,\wvec{x}}  \wfc{p_i}{t,\wvec{x},\wvec{y}}}{\sum_{i = -d}^{n} \wfc{\tilde{\lambda}_i}{t,\tilde{\wvec{x}}}} +
  \frac{ \wlr{\sum_{i = n - d}^{n} \frac{\wfc{\tilde{\lambda}_i}{t,\tilde{\wvec{x}}}}{\wfc{\tilde{\lambda}_{n-d}}{t,\wvec{x}}}}
\wfc{\tilde{\lambda}_{n-d}}{t,\wvec{x}}\wfc{p_{n-d}}{t,\wvec{x},\wvec{y}}  }{\sum_{i = -d}^{n} \wfc{\tilde{\lambda}_i}{t,\tilde{\wvec{x}}} } .
 \end{equation}
(Since $n > d + 1$ the sums in numerator of the first and last parcels in the expression above do not overlap, even
when the sum in the numerator in the middle is empty.)
Let $t \ \in \ (x_0, x_1)$ be fixed. We claim that  $y^*_0, y^*_1, \dots, y^*_{n} \ \in \ \{-1,1\}$ defined by
\[
y^*_0  := y^*_1  :=  (-1)^d \hspace{0.6cm} \wrm{and} \hspace{0.6cm}
y^*_j :=  (-1)^{d+j-1} \ \wrm{for} \ \ 2 \leq j \leq n
\]
satisfy
\begin{equation}
\label{lem2eq}
 \wfc{\tilde{\lambda}_i}{t,\wvec{x}} \wfc{p_i}{t,\wvec{x},\wvec{y^*}}  =  \sum\limits_{j = 0}^{d} \frac{\wabs{u_j}}{\wabs{t - \tilde{x}_{i+j}}}\geq 0
 \hspace{0.4cm} \wrm{for} \ i = 0 \ \mbox{or} \ 2 \leq i \leq n - d,
\end{equation}
where
\[
u_j := \frac{(-1)^{d- j}}{d! h^d}\left( \begin{array}{c} d \\ j \end{array} \right).
\]
In fact, \wref{LambdaAnd} and Equations (3.1), (3.3) and (5.1) of \cite{Ber} show that
\begin{equation}\label{And2}
 \tilde{\lambda}_i(t,\wvec{x})p_i(t,\wvec{x},\wvec{y^*})  =  (-1)^i \sum\limits_{j = 0}^{d} \frac{u_jy^*_{i+j}}{t-x_{i+j}},
\end{equation}
and \wref{lem2eq} follows from
\[
\frac{u_0y^*_0}{t-x_{0}} = \frac{\wabs{u_0}}{\wabs{t-x_{0}}}, \hspace{1cm} \frac{u_1 y^*_1}{t-x_{1}} = \frac{\wabs{u_1}}{\wabs{t-x_{1}}}
\]
and
\[
\frac{u_j y^*_{i+j}}{t-x_{i+j}} = \frac{(-1)^{d-j}\wabs{u_j}y^*_{i+j}}{t-x_{i+j}} = \frac{(-1)^{i-1}\wabs{u_j}}{t-x_{i+j}} =
\frac{(-1)^i \wabs{u_j}}{\wabs{t-x_{i+j}}}
\]
for $2 \leq i + j \leq n$ and $0 \leq j \leq d$.

Note that
\begin{equation}
\label{FHinterp4}
  \frac{\sum_{i = -d}^{0} \wfc{\tilde{\lambda}_i}{t,\tilde{\wvec{x}}}}
  {\tilde{\lambda}_0(t,\wvec{x})}  > 0,
\end{equation}
because, for $-d \leq i < 0$, we have that $d+ i+ 1 \geq 1$ and \wref{LambdaAnd} yields
\[
\frac{\wfc{\tilde{\lambda}_i}{t,\tilde{\wvec{x}}}}{\wfc{\tilde{\lambda}_0}{t,\tilde{\wvec{x}}}}
=
\frac{ \wlr{-1}^{i} \wlr{t - x_0}  \dots \wlr{t - x_{d + i}} \wlr{t - x_{d+ i + 1}} \dots \wlr{t - x_{d} }}
     { \wlr{t - \tilde{x}_i} \dots \wlr{t - \tilde{x}_{-1}} \wlr{t - x_0} \dots \wlr{t - x_{d + i}} }
\]
\begin{equation}
\label{same_sign}
= \frac{ \wlr{-1}^{i} \wlr{t - x_{d + i + 1}} \dots \wlr{t - x_{d} }}
     { \wabs{t - \tilde{x}_i} \dots \wabs{t - \tilde{x}_{-1}} }
= \frac{\wabs{t - x_{d + i + 1}} \dots \wabs{t - x_{d} }}
     { \wabs{t - \tilde{x}_i} \dots \wabs{t - \tilde{x}_{-1}} } \geq 0,
\end{equation}
and this inequality also holds for $i = 0$.
%%%%%%%%%%%%% alternate %%%%%%%%%%%%%%%%%%
Moreover, the signs of the numbers
$\tilde{\lambda}_{1}(t,\tilde{\wvec{x}})$, $\tilde{\lambda}_{2}(t,\tilde{\wvec{x}})$, $\dots$,
$\tilde{\lambda}_{n}(t,\tilde{\wvec{x}})$ alternate, and their magnitude decreases because, for $1 \leq i < n$, \wref{LambdaAnd} yields
\begin{equation}
\label{alternate}
-1 < \frac{\tilde{\lambda}_{i + 1}(t,\tilde{\wvec{x}})}{\tilde{\lambda}_{i}(t,\tilde{\wvec{x}})} = - \frac{\tilde{x}_{i} - t}{\tilde{x}_{d + i + 1} - t} < 0.
\end{equation}
As a result,
\[
 \sum_{i = n - d}^{n} \frac{\wfc{\tilde{\lambda}_i}{t,\tilde{\wvec{x}}}}{\wfc{\tilde{\lambda}_{n-d}}{t,\wvec{x}}}  \geq 0.
\]
This inequality and \wref{lem2eq} with $i = n - d$ imply that the numerator of the last parcel
in the sum in the right hand side of \wref{FHinterp3} is not negative,
and combining \wref{FHinterp3}, \wref{lem2eq} and  \wref{FHinterp4} we obtain
\begin{equation}
\label{FHinterp5}
\wabs{\wfc{\tilde{r}_{\wvec{x},d,d,d}[\wvec{y}^*]}{t}} \geq
 \frac{\frac{\sum_{i = -d}^{0} \tilde{\lambda}_i(t,\tilde{\wvec{x}})}
 {\vspace{0.5mm}\tilde{\lambda}_0(t,\wvec{x})} \tilde{\lambda}_0(t,\wvec{x})p_0(t,\wvec{x},\wvec{y^*}) -
 \wabs{\tilde{\lambda}_1(t,\wvec{x}) p_1(t,\wvec{x},\wvec{y^*})}   }
 {\wabs{\sum_{i = -d}^{n} \tilde{\lambda}_i(t,\tilde{\wvec{x}})} }.
\end{equation}
Moreover, \wref{lem2eq} and \wref{And2} yield
\[
\tilde{\lambda}_0(t,\wvec{x})p_0(t,\wvec{x},\wvec{y^*})  - \wabs{\tilde{\lambda}_1(t,\wvec{x})p_1(t,\wvec{x},\wvec{y^*})}  \ \  \geq \ \
 \sum\limits_{j = 0}^{d} \frac{\wabs{u_j}}{\wabs{t-x_{j}}} - \sum\limits_{j = 0}^{d} \frac{\wabs{u_{j}}}{\wabs{t-x_{j+1}}}
\]
\begin{equation}\label{posit}
 \geq\  \frac{\wabs{u_0}}{\wabs{t-x_0}} +
 \wset{\frac{\wabs{u_1}-\wabs{u_0}}{\wabs{t-x_1}} - \frac{\wabs{u_1}}{\wabs{t-x_2}}}
  + \sum\limits_{j = 2}^{d} \left(\frac{\wabs{u_j}}{\wabs{t-x_{j}}} -
  \frac{\wabs{u_j}}{\wabs{t-x_{j+1}}} \right).
\end{equation}
 The last sum in \wref{posit} is  positive for all $t \ \in \ (x_0, x_1)$.
 The term in brackets is also positive for $t \ \in \ (x_0, x_1)$, because
\[
 \frac{\wabs{u_1}-\wabs{u_0}}{\wabs{u_1}} = 1 - \frac{1}{d} \geq \frac{1}{2} \geq \frac{\wabs{t-x_1}}{\wabs{t-x_2}}.
\]
Therefore,
\begin{equation}
\label{pos_diff}
\tilde{\lambda}_0(t,\wvec{x})p_0(t,\wvec{x},\wvec{y^*})  - \wabs{ \tilde{\lambda}_1(t,\wvec{x})p_1(t,\wvec{x},\wvec{y^*}) }  \ \  \geq \ \ 0.
\end{equation}
Equation \wref{same_sign} shows that
the numbers $\tilde{\lambda}_{-d}(t,\tilde{\wvec{x}})$, $\tilde{\lambda}_{1-d}(t,\tilde{\wvec{x}})$, $\dots$, $\tilde{\lambda}_{-1}(t,\tilde{\wvec{x}})$
have the same sign as $\tilde{\lambda}_0(t,\tilde{\wvec{x}})$. When $-d < i \leq 0$ we also have
\[
\frac{\tilde{\lambda}_{i}(t,\tilde{\wvec{x}})}{\tilde{\lambda}_1(t,\tilde{\wvec{x}})} =
\frac{\wlr{-1}^{i + 1}\wlr{t - x_1} \dots \wlr{t - x_{d + i}} \wlr{t - x_{d + i + 1}} \dots \wlr{t - x_{d + 1}}}
{\wlr{t - \tilde{x}_i} \dots \wlr{t - x_0} \wlr{t - x_1} \dots \wlr{t - x_{d + i}}} =
\]
\begin{equation}
\label{q01}
\frac{\wlr{-1}^{i + 1}\wlr{t - x_{d + i + 1}} \dots \wlr{t - x_{d + 1}}}
{\wlr{t - \tilde{x}_i} \dots \wlr{t - x_0}}
=
\frac{\wlr{x_{d + i + 1} - t} \dots \wlr{x_{d + 1} - t}}
{\wlr{t - \tilde{x}_i} \dots \wlr{t - x_0}} \geq 0,
\end{equation}
and the reader can verify that this inequality also holds for $i = -d$. Therefore,
\begin{equation}
\label{same_signB}
\frac{\tilde{\lambda}_i(t,\tilde{\wvec{x}})}{\tilde{\lambda}_{1}(t,\tilde{\wvec{x}})} > 0 \hspace{1cm} \wrm{for} \hspace{0.2cm} -d \leq i \leq 1.
\end{equation}
Since all $\wfc{\tilde{\lambda}_i}{t,\wvec{x}}$ for $- d \leq 0 \leq 1$ have the same sign,
and for $i \geq 1$ the
signs of the $\wfc{\tilde{\lambda}_i}{t,\wvec{x}}$ alternate and their magnitude decreases, we have
\begin{equation}
\label{denominator}
\wabs{\sum_{i = -d}^n \wfc{\tilde{\lambda}_i}{t,\tilde{\wvec{x}}}}
 \leq
\wabs{\sum_{i = -d}^1 \wfc{\tilde{\lambda}_i}{t,\tilde{\wvec{x}}}}.
\end{equation}
Equations \wref{FHinterp5}, \wref{pos_diff} and \wref{denominator}  lead to
\begin{equation}
\label{bounda}
\wabs{\wfc{\tilde{r}_{\wvec{x},d,d,d}[\wvec{y}^*]}{t}} \geq
  \frac{\wabs{\sum_{i = -d}^{-1} \tilde{\lambda}_i(t,\tilde{\wvec{x}})} \wabs{p_0(t,\wvec{x},\wvec{y^*})}}{\wabs{\sum_{i = -d}^{1} \tilde{\lambda}_i(t,\tilde{\wvec{x}})}}.
\end{equation}
When $-d < i \leq 0$, Equation \wref{same_sign}, and the comment just after it, yield
\[
\frac{\wfc{\tilde{\lambda}_i}{t,\tilde{\wvec{x}}}}{\wfc{\tilde{\lambda}_0}{t,\tilde{\wvec{x}}}}
\geq
\frac{ \wlr{x_{d + i + 1} - x_1} \dots \wlr{x_d - x_1}}
     { \wlr{x_1 - \tilde{x}_i} \dots \wlr{x_1 - \tilde{x}_{-1}} }
= \frac{\wlr{d + i} \dots \wlr{d-1} h^{-i}}
     { \wlr{1 - i}! h^{-i} } =
\]
\[
\frac{\wlr{d + i} \dots \wlr{d-1}}{ \wlr{1 - i}!} = \frac{1}{d} \wvecB{d}{1 - i}.
\]
Therefore,
\begin{equation}
\label{sum0}
d \times \frac{\sum_{i = -d}^{1} \wabs{\wfc{\tilde{\lambda}_i}{t,\tilde{\wvec{x}}}}}{\wabs{\wfc{\tilde{\lambda}_0}{t,\tilde{\wvec{x}}}}}
\geq \sum_{i = 1 - d}^{0} \wvecB{d}{1 - i} = \sum_{j = 1}^{d}  \wvecB{d}{j}
= 2^d - 1.
\end{equation}
Moreover, for $-d < i \leq 0$ Equation \wref{q01} yields
\[
\frac{\tilde{\lambda}_{i}(t,\tilde{\wvec{x}})}{\tilde{\lambda}_1(t,\tilde{\wvec{x}})}
\geq
\frac{\wlr{x_{d + i + 1} - x_1} \dots \wlr{x_{d + 1} - x_1}}
{\wlr{x_1 - \tilde{x}_i} \dots \wlr{x_1 - x_0} } =
\frac{\wlr{d + i} \dots d\, h^{1 - i}}
{\wlr{1 - i}!  h^{1 - i}} = \wvecB{d}{1 - i},
\]
and
\begin{equation}
\label{sum1}
\frac{\sum_{i = -d}^1 \wabs{\tilde{\lambda}_{i}(t,\tilde{\wvec{x}})}}{\wabs{\tilde{\lambda}_1(t,\tilde{\wvec{x}})}}
\geq \sum_{i = 1- d}^0 \frac{\wabs{\tilde{\lambda}_{i}(t,\tilde{\wvec{x}})}}{\wabs{\tilde{\lambda}_1(t,\tilde{\wvec{x}})}} + 1
\geq \sum_{i = 1 - d}^1 \wvecB{d}{1 - i} =
\sum_{j = 0}^d \wvecB{d}{j} = 2^d.
\end{equation}
It follows from Equations \wref{main_thesis}, \wref{sum0} and \wref{sum1} that
\[
\wabs{\sum_{i = -d}^{-1} \tilde{\lambda}_i(t,\tilde{\wvec{x}})}
= \sum_{i = -d}^{-1} \wabs{\tilde{\lambda}_i(t,\tilde{\wvec{x}})} =
\sum_{i = -d}^{1} \wabs{\tilde{\lambda}_i(t,\tilde{\wvec{x}})} - \wabs{\tilde{\lambda}_0(t,\tilde{\wvec{x}})}  - \wabs{\tilde{\lambda}_1(t,\tilde{\wvec{x}})}
\]
\[
=
\wlr{1 - \frac{\wabs{\tilde{\lambda}_0(t,\tilde{\wvec{x}})}}{\sum_{i = -d}^{1} \wabs{\tilde{\lambda}_i(t,\tilde{\wvec{x}})} }  -
\frac{\wabs{\tilde{\lambda}_1(t,\tilde{\wvec{x}})}}{\sum_{i = -d}^{1} \wabs{\tilde{\lambda}_i(t,\tilde{\wvec{x}})} }}  \sum_{i = -d}^{1} \wabs{\tilde{\lambda}_i(t,\tilde{\wvec{x}})}
\geq
\kappa_d \wabs{\sum_{i = -d}^{1} \tilde{\lambda}_i(t,\tilde{\wvec{x}})}.
\]
The inequality in the previous line and \wref{bounda}  yield
\[
\wabs{\wfc{\tilde{r}_{\wvec{x},d,d,d}[\wvec{y}^*]}{t}} \geq  \kappa_d \wabs{p_0(t,\wvec{x},\wvec{y^*})}.
\]
Equation \wref{lem2eq} shows that,
for $t \ \in [x_0, x_1]$, $\wabs{p_0(.,\wvec{x},\wvec{y^*})}$
is identical to the Lebesgue function for polynomial interpolation at $(d+1)$
equally spaced nodes in $[x_0,x_d]$. According to \cite{Brutman}, the Lebesgue
function for polynomial interpolation at equally spaced nodes attains its maximum
at some $t^* \in (x_0,x_1)$. For this $t^*$ we have
\[
\wlebcdiag{} \geq \wlebfxfhdiag{t^*}
\ \geq \ \wabs{\wfc{\tilde{r}_{\wvec{x},d,d,d}[\wvec{y^*}]}{t^*}} \ \geq
\kappa_d \wabs{ \wfc{p_0}{t^*,\wvec{x},\wvec{y^*}}} =
   \kappa_d \Lambda_d,
\]
and this completes the proof of Theorem \ref{thm_main}. \qed{}

\section{Backward instability}
\label{section_backward}
In this section we discuss the backward stability of
the barycentric formula used to evaluate extended interpolants
when $\tilde{\wvec{y}}$ is given by a function
$\wfc{Y}{d, \wvec{x}, \wvec{y}}$ which is linear in $\wvec{y}$.
Formally, we take $\tilde{y}_i := y_i$ for $0 \leq i \leq n$ and
\begin{equation}
\label{ytilback}
\tilde{y}_i := \wfc{Y_i}{d, \wvec{x}, \wvec{y}} = \sum_{j = 0}^n \wfc{h_{ij}}{d,\wvec{x}} y_j
\hspace{1cm} \wrm{for} \hspace{0.2cm} i \in \wset{-d,\dots, n + d} \setminus \wset{0,\dots, n}.
\end{equation}
The extended function values $\tilde{\wvec{y}}$ are supposed to be evaluated numerically
and then to be used to evaluate the barycentric interpolant $\wfc{b}{t, \wvec{y}}$  given by
\begin{equation}
\label{baryback}
\wfc{b}{t, \wvec{y}}  := \left. \sum_{i = -d}^{n+d} \frac{\tilde{w}_i \tilde{y}_i}{t - \tilde{x}_i} \right/
\sum_{i = -d}^{n+d} \frac{\tilde{w}_i}{t - \tilde{x}_i}
\hspace{0.5cm} \wrm{with} \hspace{0.5cm} \tilde{y}_i = \wfc{Y_i}{d, \wvec{x}, \wvec{y}}.
\end{equation}
We assume that the weights $\tilde{w}_i$ are such that the denominator of $\wfc{b}{t,\wvec{y}}$ is different from
zero for $t \in [x_0,x_n] \setminus \wset{x_0,x_1,\dots x_n}$.

As we have shown in Section 3, the Equation (2.3) in \cite{Klein}
and the equation just before Theorem 5.1 in  \cite{BerrutKleinCAM}
are particular cases of Equation \wref{baryback}.
Therefore, by discussing the backward stability of \wref{ytilback}--\wref{baryback}
we also cover the the backward stability of the interpolation formulae proposed in the literature.

In order to analyze the backward stability of \wref{ytilback}--\wref{baryback},
it is convenient to proceed as in Section \ref{section_reduced} and rewrite \wref{baryback} as
\begin{equation}
\label{barybackred}
\wfc{b}{t, \wvec{y}} = \left. \sum_{i = 0}^{n} \wfc{d_i}{t} y_i \right/
\sum_{i = -d}^{n+d} \frac{\tilde{w}_i}{t - \tilde{x}_i},
\end{equation}
for
\begin{equation}
\label{dj}
\wfc{d_j}{t} := \sum_{i = -d}^{-1} \frac{\tilde{w}_i \wfc{h_{ij}}{d,\wvec{x}}}{t - \tilde{x}_i} + \frac{\tilde{w}_j}{t - x_j}
 + \sum_{i = n + 1}^{n + d} \frac{\tilde{w}_i \wfc{h_{ij}}{d, \wvec{x}}}{t - \tilde{x}_i}.
\end{equation}
Equations \wref{barybackred}--\wref{dj} can be verified as
in the proof of the validity of the reduced form
\wref{reduced_form} presented in Section \ref{subsection_lem_clean}.

%by analyzing the numerator $\wfc{\eta}{t,\wvec{y}}$ of the right hand
%side of \wref{baryback}
%\[
%\wfc{\eta}{t, \wvec{y}} := \sum_{i = -d}^{n+d} \frac{\tilde{w}_i \tilde{y}_i}{t - \tilde{x_i}} =
% \sum_{i = -d}^{-1} \frac{\tilde{w}_i \tilde{y}_i}{t - \tilde{x_i}} +
% \sum_{i = 0}^{n} \frac{\tilde{w}_i \tilde{y}_i}{t - \tilde{x_i}} +
% \sum_{i = n + 1}^{n+d} \frac{\tilde{w}_i \tilde{y}_i}{t - \tilde{x_i}}.
%\]
%Note that
%\[
%\wfc{\eta}{t, \wvec{y}}= \sum_{i = -d}^{-1}  \frac{\tilde{w}_i \sum_{j = 0}^n h_{ij} y_j}{t - \tilde{x_i}} +
% \sum_{i = 0}^{n} \frac{\tilde{w}_i \tilde{y}_i}{t - \tilde{x_i}} +
% \sum_{i = n + 1}^{n+d} \frac{\tilde{w}_i \sum_{j = 0}^n h_{ij} y_j}{t - \tilde{x_i}}
%\]
%\[
%= \sum_{i = -d}^{-1} \sum_{j = 0}^n \frac{\tilde{w}_i h_{ij} y_j}{t - \tilde{x_i}} +
% \sum_{j = 0}^{n} \frac{\tilde{w}_j \tilde{y}_j}{t - \tilde{x_j}} +
% \sum_{i = n + 1}^{n+d} \sum_{j = 0}^n \frac{\tilde{w}_i  h_{ij} y_j}{t - \tilde{x_i}}
%\]
%\[
%= \sum_{j = 0}^n \sum_{i = -d}^{-1}  \frac{\tilde{w}_i h_{ij} y_j}{t - \tilde{x_i}} +
% \sum_{j = 0}^{n} \frac{\tilde{w}_j y_j}{t - x_j} +
%\sum_{j = 0}^n \sum_{i = n + 1}^{n+d}  \frac{\tilde{w}_i  h_{ij} y_j}{t - \tilde{x_i}}
%\]
%\[
%= \sum_{j = 0}^n \wlr{\sum_{i = -d}^{-1}  \frac{\tilde{w}_i h_{ij}}{t - \tilde{x_i}} +
%\frac{\tilde{w}_j}{t - x_j} +
%\sum_{i = n + 1}^{n+d}  \frac{\tilde{w}_i  h_{ij}}{t - \tilde{x_i}}} y_j =
%\sum_{j = 0}^n \wfc{d_j}{t} y_j.
%\]
We adopt the definition of backward stability used by
Higham in \cite{HIGHAM_IMA}, that is, the formulae above are backward stable when the
value $\wfc{\hat{b}}{t,\wvec{y}}$ obtained by evaluating \wref{ytilback}-\wref{baryback}
in inexact arithmetic is equal to the exact value
$\wfc{b}{t,\hat{\wvec{y}}}$ for a perturbed vector $\hat{\wvec{y}}$ with
$\hat{y}_i = y_i \wlr{1 + \phi_i}$ for $\phi_i$ small.
With this in mind, we can summarize this section as follows:\\[-0.2cm]
\begin{quote}
The barycentric formula \wref{ytilback}-\wref{baryback} is not
backward stable in Higham's sense
when $\tilde{w}_k \neq 0$ for some $k \in \wset{-d,\dots,n+d} \setminus \wset{0,\dots n}$
and there exists $t \in [x_0,x_n] \setminus \wset{x_0,x_1,\dots x_n}$ and $0 \leq j \leq n$ such that
$\wfc{d_j}{t} = 0$.\\[-0.1cm]
\end{quote}

In this circumstance, we can prove the backward instability of \wref{ytilback}-\wref{baryback}
by considering $\wvec{y} \in \wrn{n+1}$ with $y_j = 1$ and $y_i = 0$ for $i \neq j$
and all $\hat{\wvec{y}} \in\wrn{n+1}$ with $\hat{y}_i = y_i \wlr{1 + \phi_i}$ for some $\phi_i \in \wrone{}$.
On the one hand, we have that $\wfc{d_i}{t} y_i \wlr{1 + \phi_i} = 0$ for all $i$ and $\phi_i \in \wrone{}$,
because $\wfc{d_j}{t} = 0$ and $y_i = 0$ for $i \neq j$.
Therefore, Equation \wref{barybackred} shows that when we evaluate \wref{ytilback}-\wref{baryback}
in exact arithmetic with $t$ and $\hat{\wvec{y}}$ we obtain
$\wfc{b}{t, \hat{\wvec{y}}} = 0$. On the other hand, if the unique rounding error occurs
in the evaluation of Equation  \wref{ytilback} for $i = k$, so that
$\tilde{y}_{k}$ is computed as  $\tilde{y}_{k} + \xi$ for $\xi \neq 0$, then
equation \wref{barybackred} shows that the
computed value $\wfc{\hat{b}}{t,\wvec{y}}$ satisfies
\[
\wfc{\hat{b}}{t, \wvec{y}} = \left. \wlr{\frac{\tilde{w}_{k} \xi }{t - \tilde{x}_{k}} + \sum_{i = -d}^{n+d} \frac{\tilde{w}_i \tilde{y}_i}{t - \tilde{x}_i} } \right/
\sum_{i = -d}^{n+d} \frac{\tilde{w}_i}{t - \tilde{x}_i}
\]
\[
= \left. \frac{\tilde{w}_{k} \xi }{t - \tilde{x}_{k}} \right/
\sum_{i = -d}^{n+d} \frac{\tilde{w}_i}{t - \tilde{x}_i}
+  \wfc{b}{t, \wvec{y}} =
\left. \frac{\tilde{w}_{k} \xi }{t - \tilde{x}_{k}} \right/
\sum_{i = -d}^{n+d} \frac{\tilde{w}_i}{t - \tilde{x}_i} \neq 0 =  \wfc{b}{t, \hat{\wvec{y}}}.
\]
Therefore, the computed value $\wfc{\hat{b}}{t, \wvec{y}}$ differs from all the exact
values $\wfc{b}{t, \hat{\wvec{y}}}$ and, according to Higham's definition,
\wref{ytilback}-\wref{baryback} is not backward stable in this case.

\begin{figure}[!h]
\includegraphics[bb=-250 0 850 700, width=7cm, height=3.2cm]{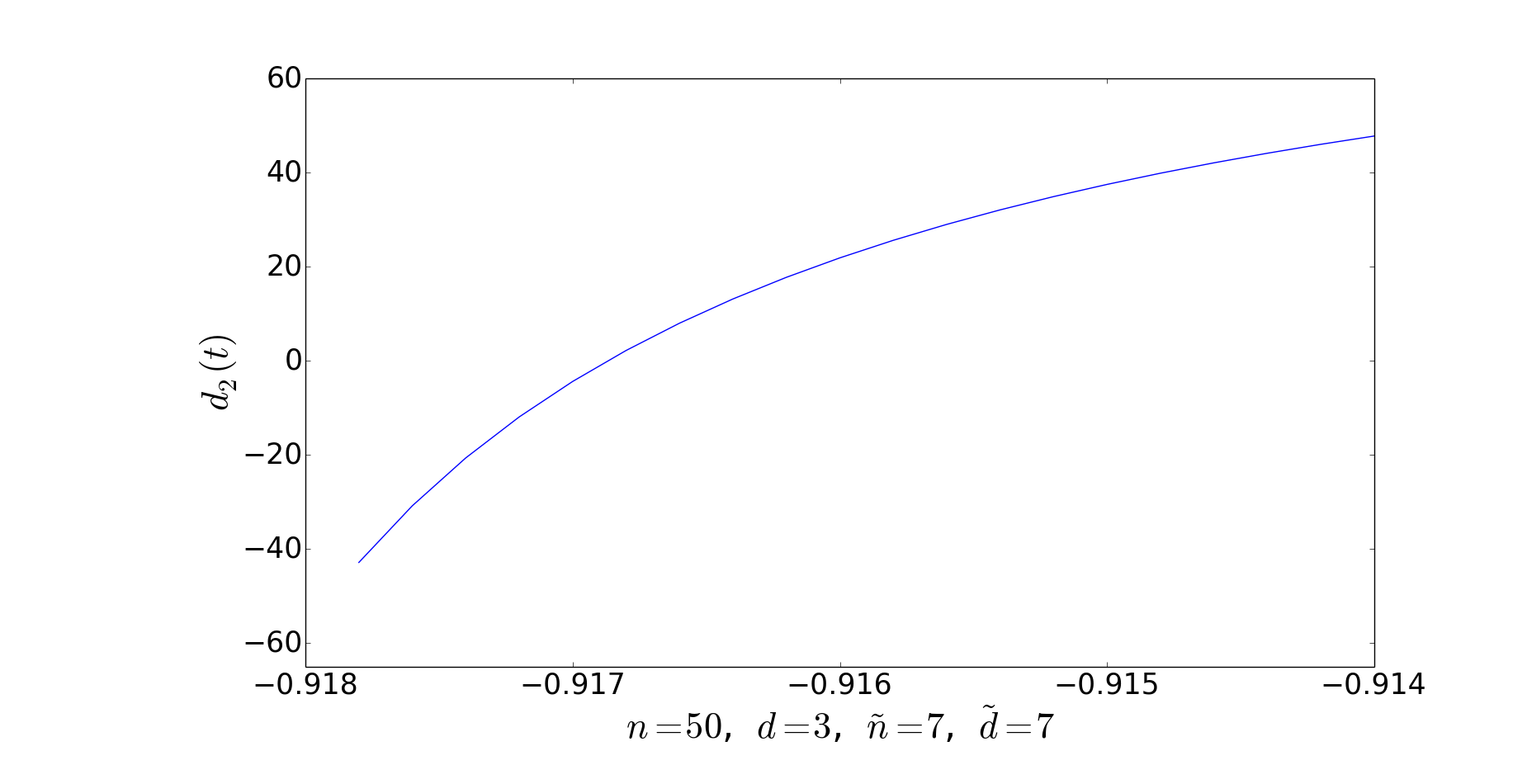}
\caption{The function $\wfc{d_2}{t}$ for the extended interpolant with $n = 50$, $d = 3$, $\tilde{n} = 11$ and $\tilde{d} = 7$
considered in \cite{Klein}.}
\label{figure_secular}
\end{figure}

In practice the $\tilde{w}_k$ are different from zero and the simple condition $\wfc{d_j}{t} = 0$ implies the backward instability
of \wref{ytilback}-\wref{baryback}. We conclude this section with
Figure \ref{figure_secular}, which shows that there is $t \in [x_0,x_n] \setminus \wset{x_0,\dots,x_n}$
for which $\wfc{d_2}{t} = 0$
for the extended interpolant with $n=50$, $d = 3$, $\tilde{n} = 11$ and
$\tilde{d} = 7$ considered by \cite{Klein}.
In fact, in this case $h = 2/50 = 0.04$, $x_{2} = -0.92$, $x_{3} = -0.88$ and
the function $\wfc{d_2}{t}$ has a zero in the interval $[-0.918,-0.914] \subset (x_2, x_3)$.

\section{Sources of numerical instability for extended interpolants}
\label{section_instability}
This section shows that extended interpolants based on extrapolation are
strongly affected by the numerical errors in the extrapolation step when $\tilde{d}$ is large.
The current literature pays little attention to this point and
presents experimental comparisons of usual and extended interpolants
that highlight cases in which $\tilde{d}$ is
much smaller than $d$. Such experiments
are biased in favor of extended interpolants: increasing
$\delta$ for usual interpolants makes as much sense as increasing
$\tilde{d}$ for extended interpolants when $d > \tilde{d}$,
because for $d > \tilde{d}$ the order of approximation of the
extended interpolants is $h^{\tilde{d} + 1}$, and not $h^{d +1}$.

\begin{figure}[!h]
\includegraphics[bb=50 0 1150 650, width=10cm, height=4.6cm]{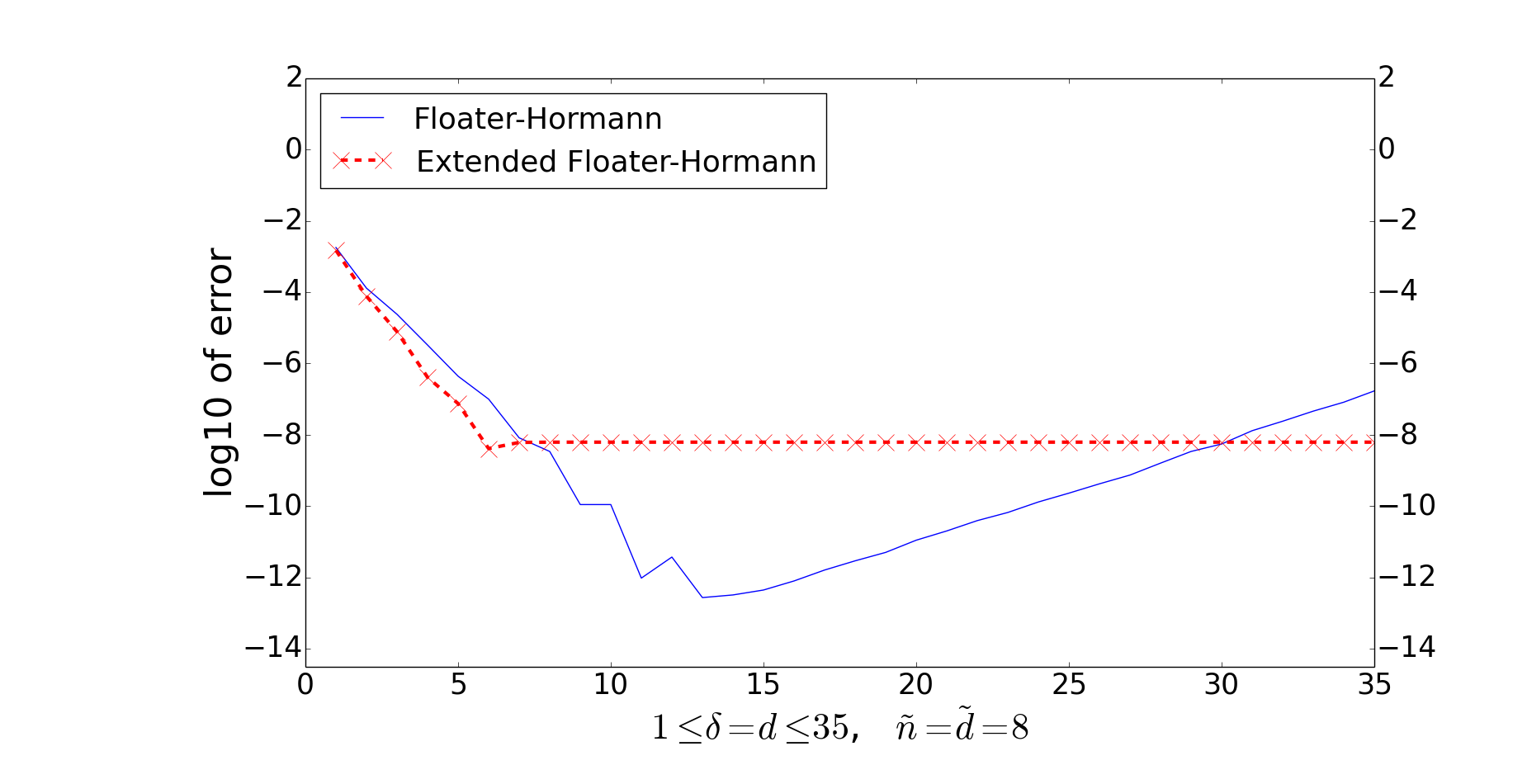}
\caption{Log10 of the error for $\wfc{f}{t} = \wfc{\sin}{20t}$ and $t \in [-1,1]$, with $n = 200$,
$\tilde{n} = \tilde{d} = 8$ and $\delta = d$ varying from $1$ to $35$. Note that by simply increasing
$d$, with an inappropriate $\tilde{d}$, we may obtain inaccurate results
for extended interpolants. In this example, increasing $d$ when $d > \tilde{d} = 8$ has no effect on the
accuracy of the extended interpolant, but increasing $\delta$ up to $13$ improves the
accuracy of the usual interpolants.
This shows that the roles of $d$ and $\delta$ are quite different when $d > \tilde{d}$.
}
\label{figure_choosing}
\end{figure}

We consider the case $\delta = d = \tilde{d} = \tilde{n}$.
In our judgment, this is the most relevant
case because it is the minimal one resulting in the same
approximation of order $h^{d+1}$ for extended interpolants and usual interpolants.
This point is reinforced by Figure \ref{figure_choosing},
which shows that it is pointless to increase $d$ when $d > \tilde{d}$.

Figure \ref{figure_choosing} illustrates the importance of choosing appropriate $\tilde{d}$
for extended interpolants and makes clear the distinction between $d$ and $\delta$.
Large values of $\delta$ have a devastating effect on usual Floater-Hormann interpolants,
and this basic fact is mentioned explicitly in the documentation of libraries that
implement these interpolants \cite{ALGLIB}. Consequently, there is
little to be learned from comparisons of extended and usual Floater-Hormann interpolants
as in Figures 6 of \cite{BerrutKleinCAM} and  \cite{Klein}: they fix $\tilde{d}$ at
small values for extended interpolants and then raise $d = \delta$ to values as large as 50.
Such choices of a large $d = \delta$ have no practical motivation for extended
interpolants and are unfavorable to usual Floater-Hormann interpolants.

\begin{figure}[!h]
\includegraphics[bb=50 0 1150 650, width=10cm, height=4.6cm]{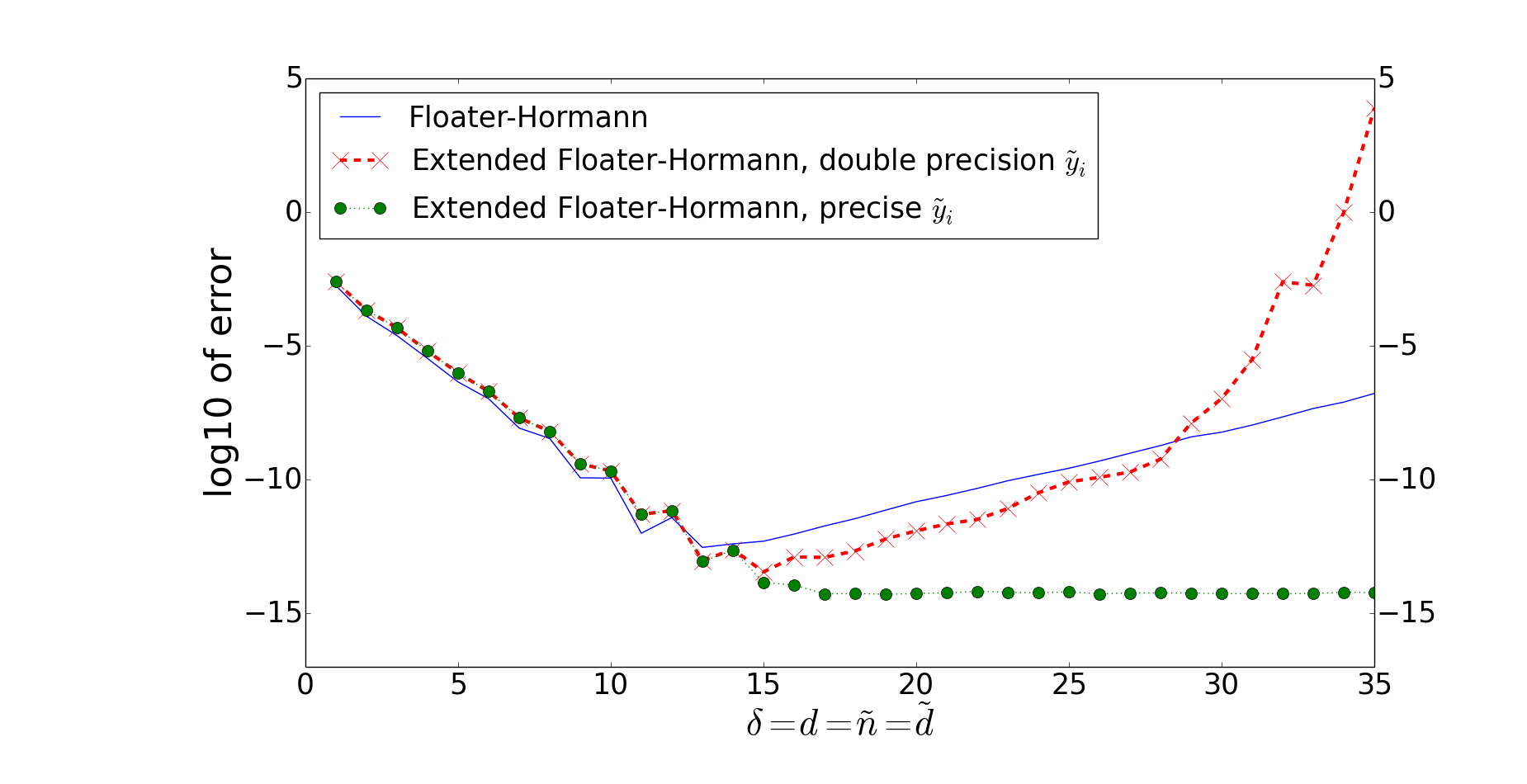}
\caption{Log10 of the error for $\wfc{f}{t} = \wfc{\sin}{20t}$, $t \in [-1,1]$ and $n = 200$.}
\label{figure_rounding}
\end{figure}

From this point to the end of this section we consider the interpolation of
$\wfc{f}{t} = \wfc{\sin}{20t}$ for $t \in [-1,1]$ with $n = 200$.
Figure \ref{figure_rounding} compares usual Floater-Hormann interpolants,
extended interpolants with $\tilde{\wvec{y}}$ computed in double precision and
extended interpolants with {\it precise} $\tilde{\wvec{y}}$. By {\it precise} we mean
that $\tilde{\wvec{y}}$ was computed using the MPFR library \cite{MPFR},
with floating point numbers with a mantissa of 320 bits, from $\wvec{y}$ computed
with the same high precision.

By {\it error} in our plots we mean
the maximum difference between the numerically evaluated interpolant and the original function
at $10^7$ equally spaced points in $[-1,1]$.
The barycentric formula \wref{bary_fh} and \wref{bary_xfh} were
evaluated in double precision
($\epsilon \approx 10^{-16}$), using straightforward C++ code. The $\wvec{y}$ and $\tilde{\wvec{y}}$
computed in multiple precision were rounded to double precision and the barycentric formula corresponding to
them was also evaluated in double precision. In other words, the case {\it precise $\tilde{\wvec{y}}$}
differs from the other cases only by the precision of the $\tilde{\wvec{y}}$, and not by the precision
used to evaluate the barycentric formulae \wref{bary_fh} and \wref{bary_xfh}.

Figure \ref{figure_rounding} shows that, when the $\tilde{y}_i$ are evaluated in double
precision, extended interpolants are not significantly more accurate than usual ones
with $\delta= d$, and they become more unstable as $d = \tilde{d} = \tilde{n}$ grows.
By contrast,
extended interpolants with precise $\tilde{\wvec{y}}$ are remarkably accurate,
even for large values of $d = \tilde{d} = \tilde{n}$. This suggests that the inaccuracy
of $\tilde{\wvec{y}}$ is the cause of the numerical instability of extended interpolants
for large $d = \tilde{d} = \tilde{n}$.

\begin{figure}[!h]
\includegraphics[bb=50 0 1150 650, width=10cm, height=4.6cm]{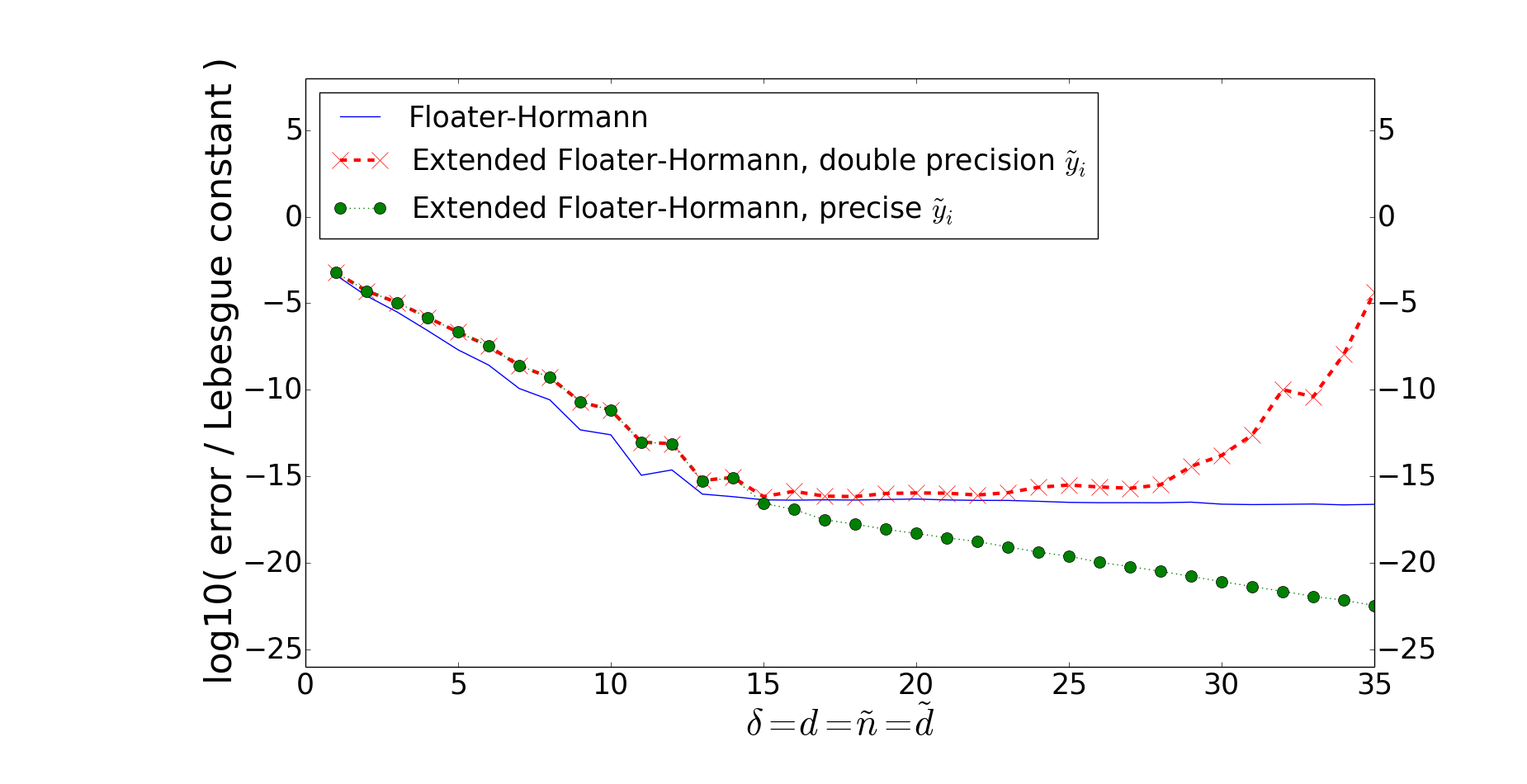}
\caption{Log10 of (forward error divided by the Lebesgue constant) for
$\wfc{f}{t} = \wfc{\sin}{20t}$, $t \in [-1,1]$ and $n = 200$.}
\label{figure_rounding_leb}
\end{figure}

Figure \ref{figure_rounding_leb} considers the ratio ``error divided by Lebesgue constant.''
This ratio is relevant for the understanding of the backward stability of interpolation formulae.
As we explain in \cite{Andre}, it is possible to implement the usual
Floater-Hormann interpolants so that the backward error is of order $n \epsilon$. Backward
errors of this order lead to forward errors of order $n \epsilon \Lambda_{\wvec{x},\delta} \wnorm{\wvec{y}}_\infty$,
 as one can verify by looking at Figures \ref{figure_rounding_leb}
and \ref{figure_rounding_semi_leb} (in this article we refer to the forward error simply
as error.) Therefore, in this case by dividing the error by the Lebesgue
constant we obtain an estimate of the backward error.
Unfortunately, Figure \ref{figure_rounding_leb} shows that the relation
between rounding errors and the Lebesgue constant for large values of $d = \tilde{d} = \tilde{n}$
for extended interpolant is more complex than the analogous relation for usual
Floater-Hormann interpolants with $\delta = d$. As a result, the fact that extended interpolants have
a smaller Lebesgue constant does not imply that they are more stable for
large values of $\delta = d = \tilde{d} = \tilde{n}$ (see Figure 5.) In fact,
in this scenario the effects of the large Lebesgue constants are quite different for extended and usual
Floater-Hormann interpolants.

The combination of Figures \ref{figure_lebesgue_constants} and \ref{figure_rounding}
leads to a surprising observation: according to Figure \ref{figure_lebesgue_constants},
the Lebesgue constant of the usual interpolants in our experiments is about $100$ times larger
than the Lebesgue constant of the corresponding extended interpolants for large
$\delta = d = \tilde{d} = \tilde{n}$, yet Figure \ref{figure_rounding} shows that
the usual interpolants are much more accurate for such large $\delta$, $d$, $\tilde{n}$ and $\tilde{d}$.

The data in Figures \ref{figure_lebesgue_constants} and \ref{figure_rounding} are combined
in Figure \ref{figure_rounding_leb}, which highlights important points for the case $d = \tilde{n} = \tilde{d} \geq 15$.
Note that this case is covered by the theory in \cite{BerrutKleinCAM} and \cite{Klein} and
their experiments consider $0 \leq d \leq 50$.
Moreover, extended interpolants are claimed to be better than usual ones for allowing
larger values of $d$ and, in view of Figure \ref{figure_choosing}, the use of
a larger $\tilde{d} \leq \tilde{n}$  is natural in this context. According to Figure \ref{figure_rounding_leb},
\begin{enumerate}
\item For fixed $n$, the error incurred by usual Floater-Hormann interpolants is of order $n \epsilon \wlebcfh{} \wnorm{f}_\infty$.\\[-0.15cm]
\item For fixed $n$, the error incurred by extended interpolants grows faster than $n \epsilon \wlebcxfh{}  \wnorm{f}_\infty$
when $d = \tilde{d} = \tilde{n}$.\\[-0.15cm]
\item The effects of the large Lebesgue constant are reduced when $\tilde{\wvec{y}}$
is precise. In this case, numerical errors occur mostly in the evaluation of the barycentric formula, and
the argument following Equation \kleinsbug{} in \cite{Klein} applies and explains
the errors of order $\epsilon$ for the extended interpolants
with precise $\tilde{\wvec{y}}$ in Figure \ref{figure_rounding}.
\end{enumerate}

\begin{figure}[!h]
\includegraphics[bb=50 0 1150 700, width=10cm, height=3.6cm]{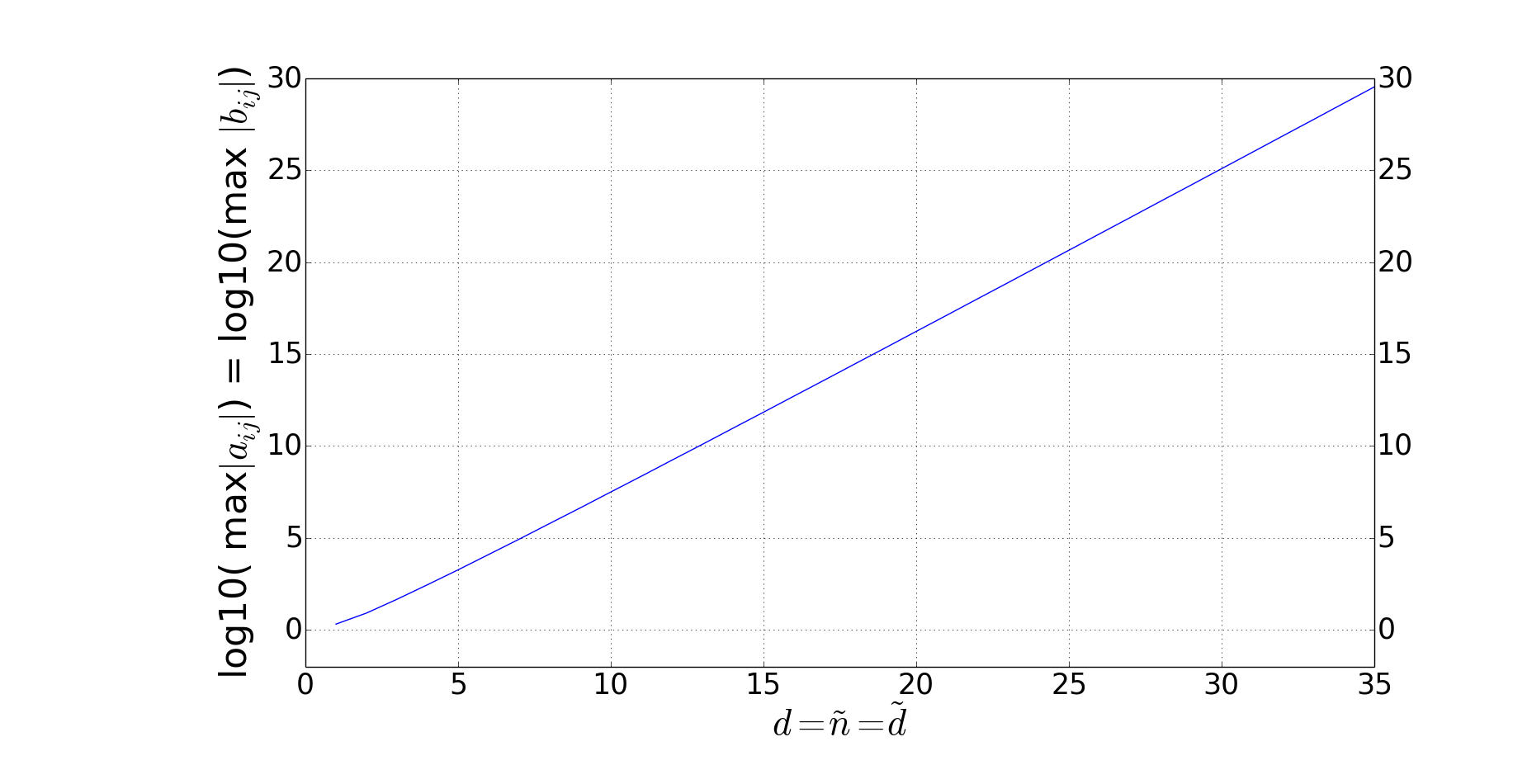}
\caption{Log10($\max\wabs{a_{ij}}$) = Log10($\max\wabs{b_{ij}}$) for $n = 200$.}
\label{figure_norm_aij}
\end{figure}

The data for the extended interpolants with double precision $\tilde{y}_i$ in Figure
\ref{figure_rounding_leb} for $d \geq 30$ indicate the existence of another source of
numerical instability for these interpolants, in addition to the
large Lebesgue constants. This extra source of instability
are the enormous entries of the matrices $a_{ij}$ and $b_{ij}$ in
Lemma \ref{lem_clean}, which are defined explicitly in \wref{def_aij}--\wref{def_bij}.
In fact, Figure \ref{figure_norm_aij} shows that the $a_{ij}$ and $b_{ij}$
grow exponentially with $d = \tilde{d} = \tilde{n}$.

\begin{figure}[!h]
\includegraphics[bb=50 0 1150 650, width=10cm, height=4.6cm]{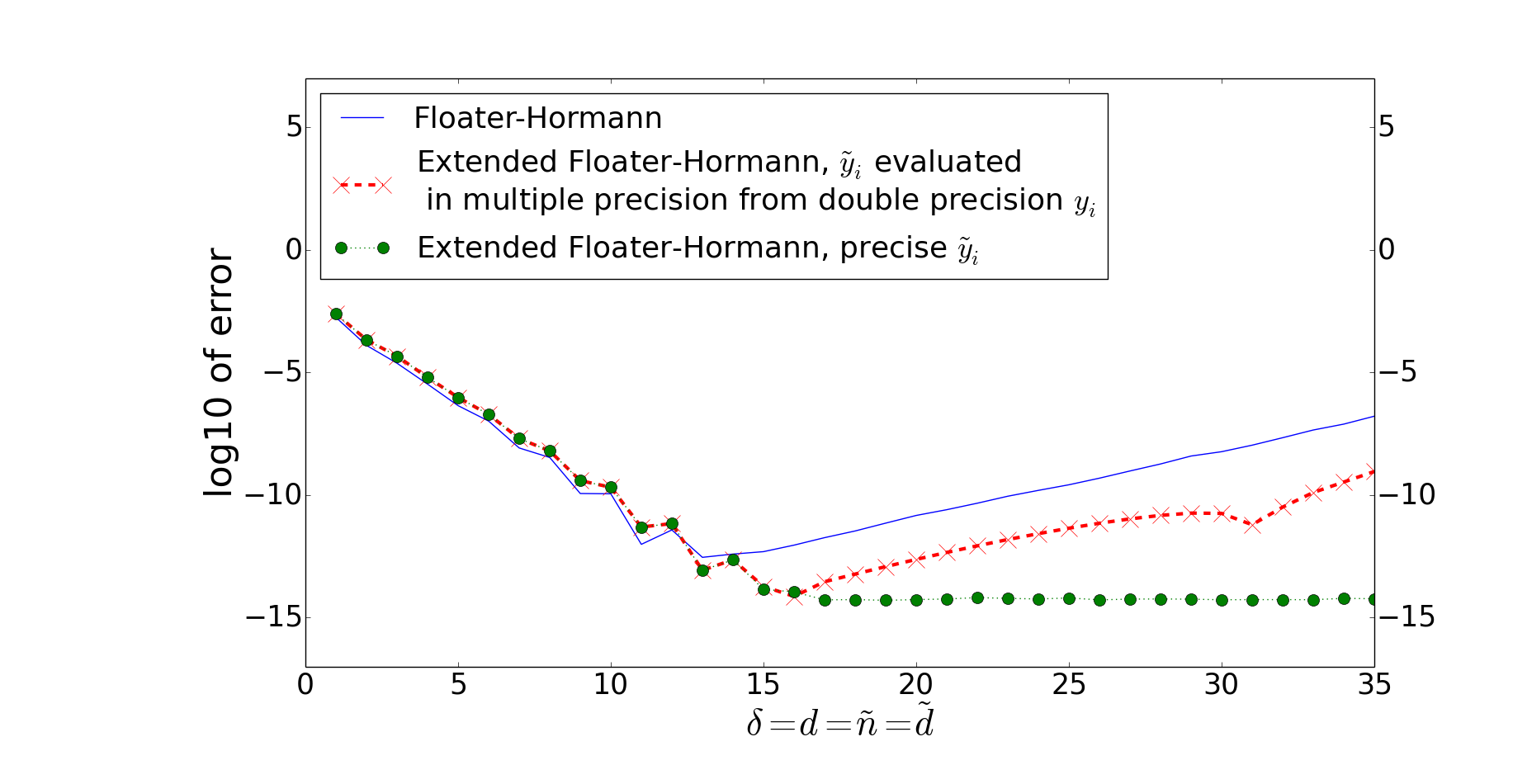}
\caption{Log10 of the forward error for $\wfc{f}{t} = \wfc{\sin}{20t}$,  $ t \in [-1,1]$ and $n = 200$.}
\label{figure_roundingsemi}
\end{figure}

\begin{figure}[!h]
\includegraphics[bb=50 0 1150 650, width=10cm, height=4.6cm]{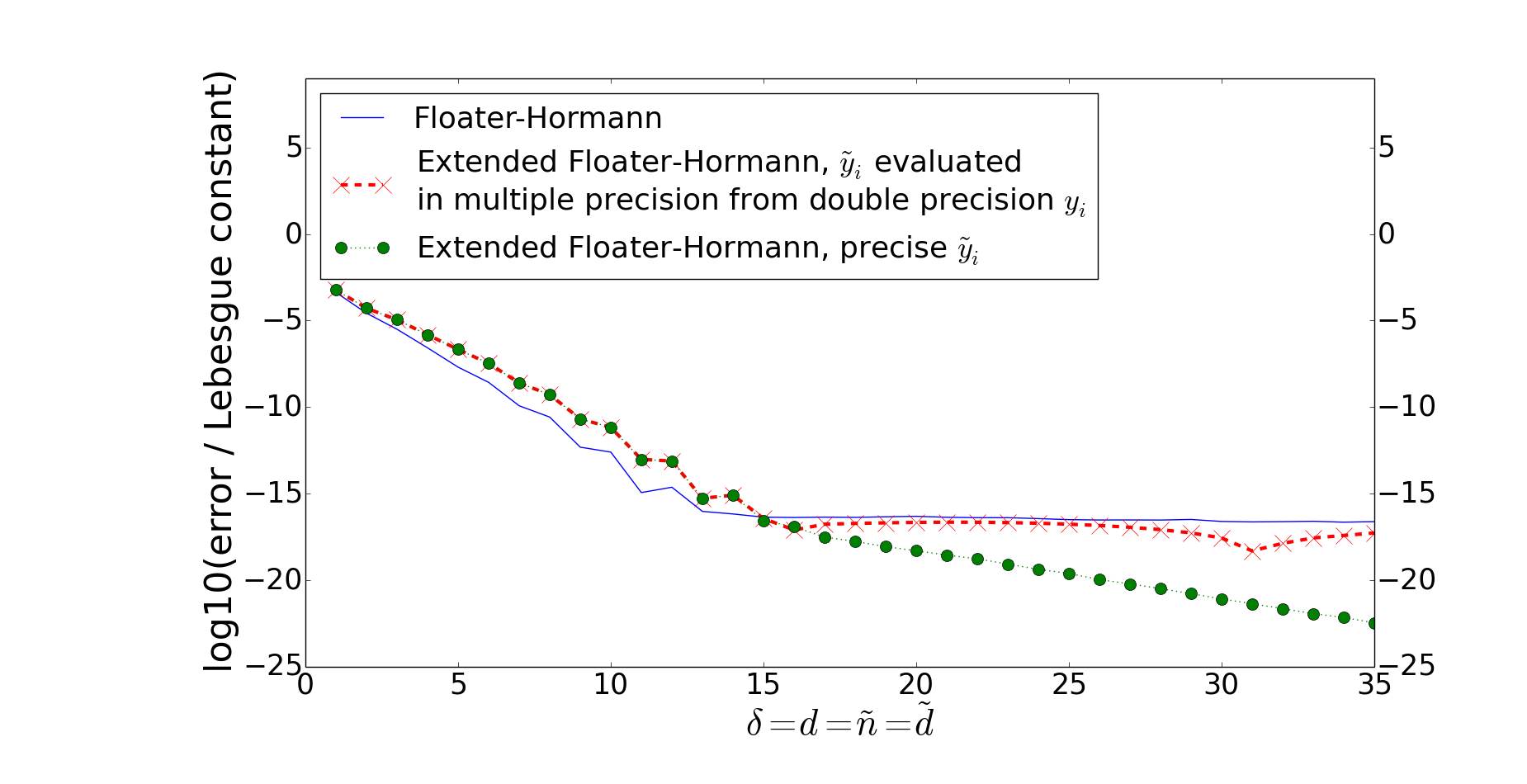}
\caption{Log10 of (forward error divided by the Lebesgue constant) for $\wfc{f}{t} = \wfc{\sin}{20t}$,
$t \in [-1,1]$ and $n = 200$.}
\label{figure_rounding_semi_leb}
\end{figure}

In view of the remarkable accuracy of the extended interpolants with precise $\tilde{\wvec{y}}$,
it  makes sense to consider the possibility of using $\tilde{\wvec{y}}$
evaluated in multiple precision from the double precision $y_i$ which are usually
available. These $\tilde{\wvec{y}}$ are not as precise as the ones obtained from
high precision $\wvec{y}$ using multiple precision arithmetic.
Figure \ref{figure_roundingsemi} illustrates, however, that this strategy
improves the accuracy of
extended interpolants, at a relatively low cost when $d$, $\tilde{d}$ and $\tilde{n}$ are
small compared to $n$ and we want to evaluate the interpolants for many values of $t$.

In the case considered in this section,
Figure \ref{figure_rounding_semi_leb} shows that
$\tilde{\wvec{y}}$ evaluated using multiple precision, from double precision $\wvec{y}$,
lead to overall numerical errors of order $n \epsilon \wlebcxfh{}$, which
are  about 100 times smaller than the errors incurred by the usual Floater-Hormann
interpolants in our experiments for large $\delta = d = \tilde{n} = \tilde{d}$.

\appendix
\section{What can we prove about the stability of extended interpolants}
\label{section_cancellation}

\begin{figure}[!h]
\includegraphics[bb=50 0 1150 650, width=10cm, height=4.6cm]{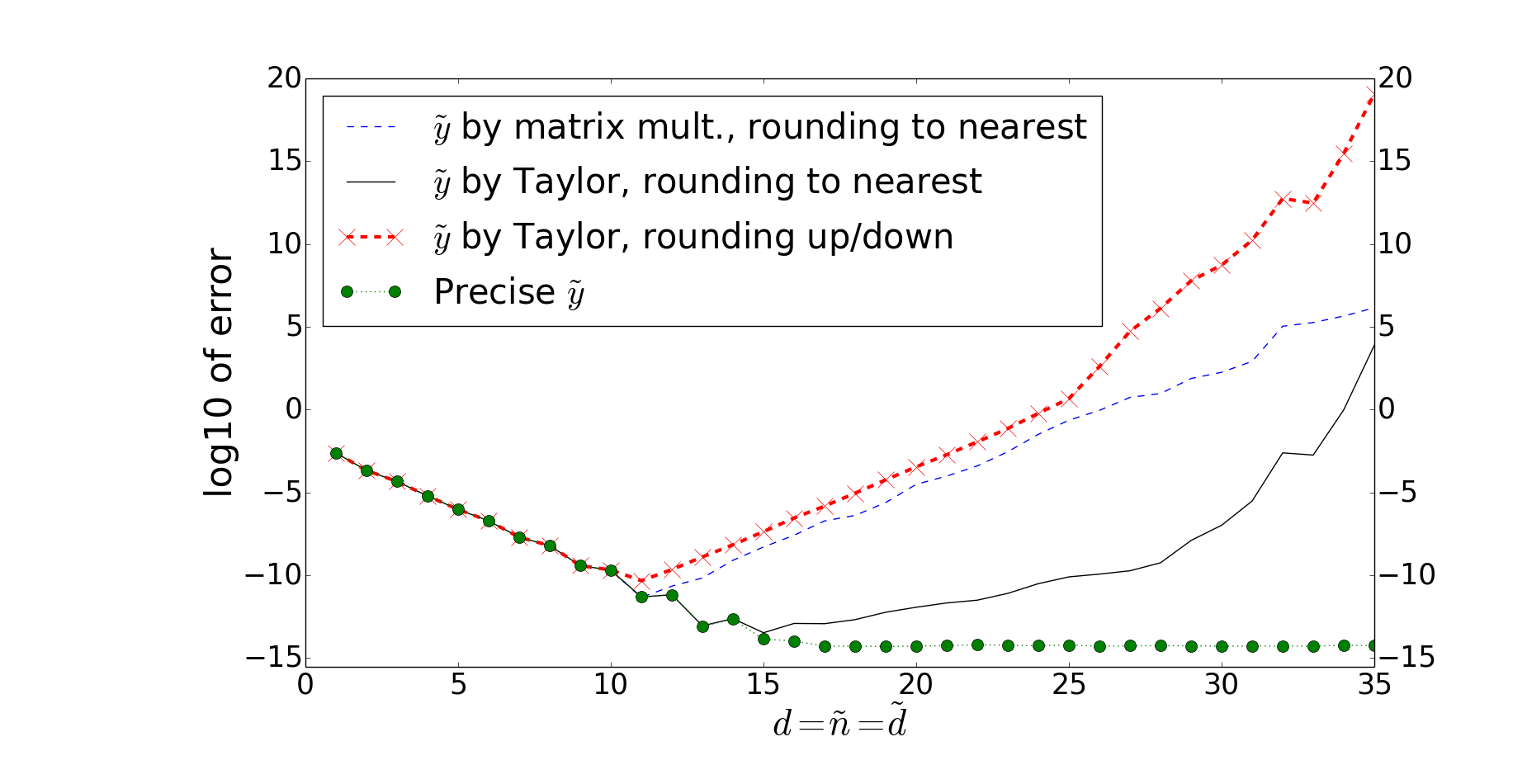}
\caption{Extended Floater-Hormann. Log10 of the forward error for $\wfc{f}{t} = \wfc{\sin}{20t}$ for
$t \in [-1,1]$ and $n = 200$. By ``$\tilde{y}$ by Taylor'' we mean $\tilde{y}$ computed
by Taylor series as in \cite{BerrutKleinCAM} and \wref{ytila}--\wref{ytilc}, and by
``$\tilde{y}$ by matrix mult.'' we mean $\tilde{y}$ computed by matrix multiplication
of $\wvec{y}$ by the matrices with entries $a_{ij}$ and $b_{ij}$, as in \wref{tilyi_clean_a}--\wref{tilyi_clean_b}.}
\label{figure_rounding_up}
\end{figure}

\begin{figure}[!h]
\includegraphics[bb=50 0 1150 650, width=10cm, height=4.6cm]{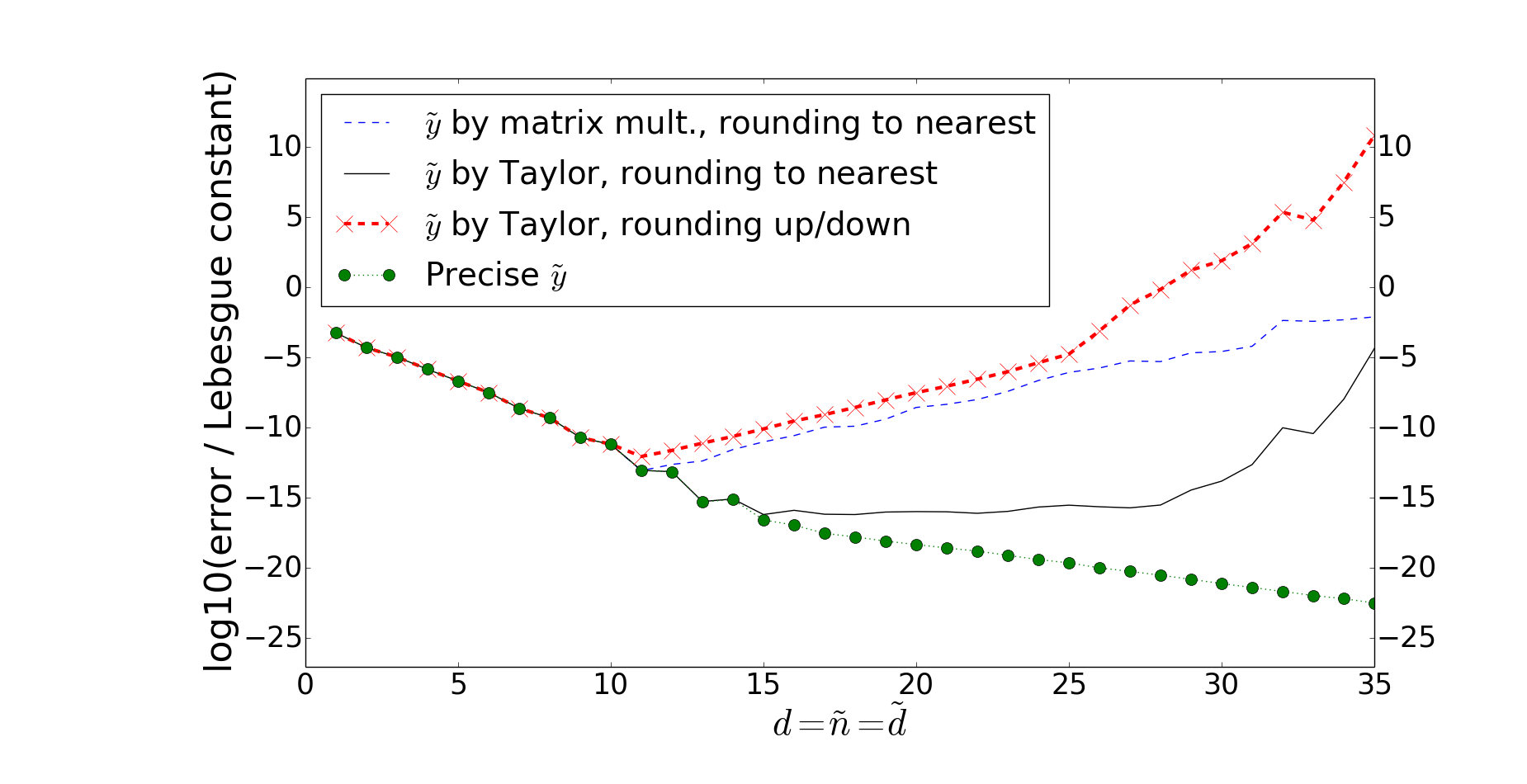}
\caption{Extended Floater-Hormann. Log10 of (forward error divided by the Lebesgue constant) for $\wfc{f}{t} = \wfc{\sin}{20t}$,
$t \in [-1,1]$ and $n = 200$.}
\label{figure_rounding_up_leb}
\end{figure}

This appendix illustrates the difficulties in building a general and realistic
stability theory for extended interpolants, a theory which would take
into account the errors introduced by the current implementations of
floating point arithmetic. We explain that the stability of extended interpolants
is sensitive to the way we implement the extrapolation step, and that
the accuracy of this step depends on the cancellation of the rounding errors.
In fact, the errors incurred by extended interpolants
can be enormous when we use the extrapolation formula proposed in
\cite{BerrutKleinCAM} and \cite{Klein},
and compute $\tilde{\wvec{y}}$ according to the following procedure:
\begin{itemize}
\item[(a)] If $i$ is even, set the rounding mode upward and evaluate
$\tilde{y}_i$ as in \wref{ytila}--\wref{ytilc}.
\item[(b)] If $i$ is odd, set the rounding mode downward and evaluate
$\tilde{y}_i$ as in \wref{ytila}--\wref{ytilc}.
\end{itemize}
In this scenario the overall effect of rounding errors can be much larger
than what one would expect from the already large Lebesgue constants,
as illustrated in Figures \ref{figure_rounding_up} and \ref{figure_rounding_up_leb}.
In the plots corresponding to $\tilde{\wvec{y}}$
evaluated as in \wref{tilyi_clean_a}--\wref{tilyi_clean_b} in these figures,
 $\tilde{\wvec{y}}$ was obtained by matrix multiplication,  with
$a_{ij}$ and $b_{ij}$ computed in multiple precision and then rounded to double precision
i.e., with $a_{ij}$ and $b_{ij}$ as accurate as possible.

We emphasize that the choices of rounding modes in the steps (a) and (b) above are not frivolous. Their purpose
is to help us understand what can be proved about the numerical stability of extended interpolants,
so that we do not try to prove something that cannot be proved.
It is  unlikely that $\tilde{\wvec{y}}$ will be evaluated as in the steps (a) and (b) when rounding to nearest.
In this mode, there is a 50\% chance of rounding up in each flop and,
under the questionable hypothesis of independence of the rounding errors,
there would be a minuscule probability of $2^{-2\wlr{d+1} d}$ of having all the intermediate
results in the evaluation of $\tilde{\wvec{y}}$ rounded up when rounding to nearest.
Since the set of floating point numbers is finite, such a coincidence may be impossible.
However, our experiments indicate that it is difficult to build a realistic theory on
the effects of rounding errors on extended interpolants, because the rounding errors induced
by our changes of rounding modes would be allowed by usual
models of floating point arithmetic, with $\epsilon$ replaced by $2 \epsilon$.
More precisely: when evaluating $\tilde{\wvec{y}}$, with our choices of rounding modes,
we monitored the relative errors
\[
\wabs{\frac{\wrounde{x + y} - \wlr{x + y}}{\wlr{x + y}}} \hspace{1cm} \wrm{and} \hspace{1cm}
\wabs{\frac{\wrounde{x * y} - \wlr{x * y}}{\wlr{x * y}}}
\]
for each operation we performed, and found all of them to be smaller than $1.97 \epsilon$.

In other words, a stability theory for extended interpolants
based on the usual models of floating point arithmetic would need
to cover the changes of rounding modes in steps (a) and (b) above
and, as a result, its predictions would be too pessimistic. Therefore, a realistic stability
theory for extended interpolants will require additional hypothesis regarding the floating
point arithmetic. By contrast, under the usual models of floating point arithmetic \cite{HIGHAM}, we
already have realistic theories bounding the rounding errors in terms of $\epsilon$, $n$ and the Lebesgue
constant for other barycentric interpolation schemes, as in \cite{Andre}
\cite{HIGHAM_IMA} \cite{Masc} \cite{MascCam} \cite{MascCamB}.

\bibliographystyle{amsplain}
\end{document}